\documentclass[12pt]{article}  
\usepackage{amsthm,amsmath,amsfonts,epsfig,amssymb,latexsym}
\usepackage{mathrsfs}
\usepackage[T1]{fontenc}
\usepackage{graphicx}
\usepackage{enumerate}
\usepackage{xy,xypic}
\usepackage{graphics}
\usepackage{footnote}
\usepackage{scalefnt}
\usepackage{tikz}
\usetikzlibrary{shapes,arrows,fit,calc,positioning}
\usetikzlibrary{matrix}
\xyoption{all}
\title{{\Large The Monoid Structure on Homotopy Obstructions} } 
\author{
 Satya Mandal\footnote{Partially supported by a General Research Grant (no 2301857) from U. of Kansas}
 ~~~and ~~Bibekananda Mishra
 \\ 
{\small University of Kansas, Lawrence, Kansas 66045, USA}\\
{\small {\it  mandal@ku.edu, bibekanandamishra@ku.edu} 
  }\\
 } 
\begin{document}
\renewcommand{\baselinestretch}{1.255}
\setlength{\parskip}{1ex plus0.5ex}
\date{23 February 2019}
\newtheorem{theorem}{Theorem}[section]
\newtheorem{proposition}[theorem]{Proposition}
\newtheorem{lemma}[theorem]{Lemma}
\newtheorem{corollary}[theorem]{Corollary}
\newtheorem{construction}[theorem]{Construction}
\newtheorem{notations}[theorem]{Notations}
\newtheorem{question}[theorem]{Question}
\newtheorem{example}[theorem]{Example}
\newtheorem{definition}[theorem]{Definition} 
\newtheorem{conjecture}[theorem]{Conjecture} 
\newtheorem{remark}[theorem]{Remark} 
\newtheorem{statement}[theorem]{Statement}

\newcommand{\iso}{\stackrel{\sim}{\longrightarrow}}

\newcommand{\sur}{\twoheadrightarrow}
\newcommand{\bD}{\begin{definition}}
\newcommand{\eD}{\end{definition}}
\newcommand{\bP}{\begin{proposition}}
\newcommand{\eP}{\end{proposition}}
\newcommand{\bL}{\begin{lemma}}
\newcommand{\eL}{\end{lemma}}
\newcommand{\bT}{\begin{theorem}}
\newcommand{\eT}{\end{theorem}}
\newcommand{\bC}{\begin{corollary}}
\newcommand{\eC}{\end{corollary}} 
\newcommand{\eop}{\hfill \rule{2mm}{2mm}}
\newcommand{\pf}{\noindent{\bf Proof.~}}
\newcommand{\PD}{\text{proj} \dim}
\newcommand{\lra}{\longrightarrow}
\newcommand{\hra}{\hookrightarrow}
\newcommand{\llra}{\longleftrightarrow}
\newcommand{\Lra}{\Longrightarrow}
\newcommand{\Llra}{\Longleftrightarrow}
\newcommand{\bE}{\begin{enumerate}}
\newcommand{\eE}{\end{enumerate}}
\newcommand{\Sets}{\underline{{\mathrm Sets}}}
\newcommand{\Sch}{\underline{{\mathrm Sch}}}
\newcommand{\ForMe}{\noindent\TCP{{\bf Remarks To Be Removed:~}}}
\newcommand{\pic}{The proof is complete.}
\newcommand{\tcp}{This completes the proof.}

\def\spec#1{\mathrm{Spec}\left(#1\right)}
\def\m{\mathfrak {m}}
\def\CA{\mathcal {A}}
\def\CB{\mathcal {B}}
\def\CP{\mathcal {P}}
\def\CC{\mathcal {C}}
\def\CD{\mathcal {D}}
\def\CE{\mathcal {E}}
\def\CF{\mathcal {F}}
\def\CE{\mathcal {E}}
\def\CG{\mathcal {G}}
\def\CH{\mathcal {H}}
\def\CI{\mathcal {I}}
\def\CJ{\mathcal {J}}
\def\CK{\mathcal {K}}
\def\CL{\mathcal {L}}
\def\CM{\mathcal {M}}
\def\CN{\mathcal {N}}
\def\CO{\mathcal {O}}
\def\CP{\mathcal {P}}
\def\CQ{\mathcal {Q}}
\def\CR{\mathcal {R}}
\def\CS{\mathcal {S}}
\def\CT{\mathcal {T}}
\def\CU{\mathcal {U}}
\def\CV{\mathcal {V}}
\def\CW{\mathcal {W}}
\def\CX{\mathcal {X}}
\def\CY{\mathcal {Y}}
\def\CZ{\mathcal {Z}}

\newcommand{\smallcirc}[1]{\scalebox{#1}{$\circ$}}
\def\BA{\mathbb {A}}
\def\BB{\mathbb {B}}
\def\BC{\mathbb {C}}
\def\BD{\mathbb {D}}
\def\BE{\mathbb {E}}
\def\BF{\mathbb {F}}
\def\BG{\mathbb {G}}
\def\BH{\mathbb {H}}
\def\BI{\mathbb {I}}
\def\BJ{\mathbb {J}}
\def\BK{\mathbb {K}}
\def\BL{\mathbb {L}}
\def\BM{\mathbb {M}}
\def\BN{\mathbb {N}}
\def\BO{\mathbb {O}}
\def\BP{\mathbb {P}}
\def\BQ{\mathbb {Q}}
\def\BR{\mathbb {R}}
\def\BS{\mathbb {S}}
\def\BT{\mathbb {T}}
\def\BU{\mathbb {U}}
\def\BV{\mathbb {V}}
\def\BW{\mathbb {W}}
\def\BX{\mathbb {X}}
\def\BY{\mathbb {Y}}
\def\BZ{\mathbb {Z}}

\newcommand{\TCP}{\textcolor{purple}}
\newcommand{\TCM}{\textcolor{magenta}}
\newcommand{\TCR}{\textcolor{red}}
\newcommand{\TCB}{\textcolor{blue}}
\newcommand{\TCG}{\textcolor{green}}

\def\SA{\mathscr {A}}
\def\SB{\mathscr {B}}
\def\SC{\mathscr {C}}
\def\SD{\mathscr {D}}
\def\SE{\mathscr {E}}
\def\SF{\mathscr {F}}
\def\SG{\mathscr {G}}
\def\SH{\mathscr {H}}
\def\SI{\mathscr {I}}
\def\SJ{\mathscr {J}}
\def\SK{\mathscr {K}}
\def\SL{\mathscr {L}}
\def\SN{\mathscr {N}}
\def\SO{\mathscr {O}}
\def\SP{\mathscr {P}}
\def\SQ{\mathscr {Q}}
\def\SR{\mathscr {R}}
\def\SS{\mathscr {S}}
\def\ST{\mathscr {T}}
\def\SU{\mathscr {U}}
\def\SV{\mathscr {V}}
\def\SW{\mathscr {W}}
\def\SX{\mathscr {X}}
\def\SY{\mathscr {Y}}
\def\SZ{\mathscr {Z}}

\def\bfA{{\bf A}}
\def\bfB{{\bf B}} 
\def\bfC{{\bf C}} 
\def\bfD{{\bf D}} 
\def\bfE{{\bf E}} 
\def\bfF{{\bf F}} 
\def\bfG{{\bf G}} 
\def\bfH{{\bf H}} 
\def\bfI{{\bf I}} 
\def\bfJ{{\bf J}} 
\def\bfK{{\bf K}} 
\def\bfL{{\bf L}} 
\def\bfM{{\bf M}} 
\def\bfN{{\bf N}} 
\def\bfO{{\bf O}} 
\def\bfP{{\bf P}} 
\def\bfQ{{\bf Q}} 
\def\bfR{{\bf R}} 
\def\bfS{{\bf S}} 
\def\bfT{{\bf T}} 
\def\bfU{{\bf U}} 
\def\bfV{{\bf V}} 
\def\bfW{{\bf W}} 
\def\bfX{{\bf X}} 
\def\bfY{{\bf Y}} 
\def\bfZ{{\bf Z}} 

\maketitle
\noindent{\bf Abstract:~}{\it 
Let $A$ be a commutative noetherian ring, containing a field $k$, with $1/2\in k$, $\dim A=d$, and let $P$ be a projective $A$-module, with $rank(P)=n$. 
Let ${\mathcal LO}(P)$ denote the set of all pairs $(I, \omega)$, where $I$ is an ideal of $A$ and $\omega: P\sur I/I^2$ is a surjective map. 
The homotopy relations on ${\mathcal LO}(P)$, induced by ${\mathcal LO}(P[T])$, leads to a set 
$\pi_0\left({\mathcal LO}(P)\right)$ of equivalence classes in ${\mathcal LO}(P)$. There are two distinguished elements 
${\bf e}_0, {\bf e}_1\in \pi_0\left({\mathcal LO}(P)\right)$, respectively, the images of $(0, 0)$ and $(A, 0)$. Define the obstruction class 
$$
\varepsilon(P)={\bf e}_0\in \pi_0\left({\mathcal LO}(P)\right),
$$
to be called the {\bf (Nori) homotopy class} of $P$.
The following results are under suitable smoothness or regularity hypotheses.
 We prove, if $2n\geq d+2$, then 
 $\pi_0\left({\mathcal LO}(P)\right)$ has a natural  structure of a monoid, which is a group if $P\cong Q\oplus A$. 
When
 $2n\geq d+3$, we prove 
 $$
 P\cong Q\oplus A \Llra \varepsilon(P)={\bf e}_1 ~~{\rm (" the~ additive ~zero")}.
 $$
Further, we give a definition of a Euler class group $E(P)$. Under suitable smoothness hypotheses, we prove, if $P\cong Q\oplus A$ and $2n\geq d+3$,
then there is natural
isomorphism $E(P) \iso \pi_0\left({\mathcal LO}(P)\right)$ of groups.
}

\newpage
\section{Introduction}

Throughout this article, unless further qualifications are added, $A$ will denote a noetherian commutative ring, with 
$\dim A=d$. Also, $P$ will denote a projective $A$-module with $rank(P)=n$.

This article is a continuation of the study of Homotopy obstructions of projective modules, that was started in \cite{MM1}.
It was pointed out in \cite{MM1} that, the study of the Homotopy obstructions of projective modules, evolved out of some germs of ideas, in two components,
 given  by Madhav V. Nori 
(around 1990), 
through some  verbal and 
 informal communications, and was referred to as the "Homotopy Program".
%
The readers would be very well advised to familiarize themselves with the introduction of [MM].  We would try to avoid any repetition,
and pick up from where we left in \cite{MM1}. We  make additional introductory comments here  only to reestablish the context. One of the two components of
these germs,
was the Homotopy Question. 
The following is a statement of the
same from \cite{M2},
which would almost certainly be an adaptation by the respective author \cite{M2}, of the more precise formulation  communicated by Nori.

\begin{question}
[Homotopy Question]\label{homoConj}
Suppose $X=\spec{A}$ is a smooth affine variety, with $\dim X=d$. Let $P$ be a projective $A$-module of rank $n$ and $f_0:P\sur I$ 
be  a surjective 
homomorphism, onto an ideal $I$ of $A$. Assume $Y=V(I)$ is smooth with $\dim Y=d-n$.
Also suppose $Z=V(J)\subseteq \spec{A[T]}=X\times \BA^1$ is a smooth subscheme, such that $Z$ intersects $X\times 0$ transversally in $Y\times 0$.
Now, suppose that $\varphi: P[T] \sur \frac{J}{J^2}$ is a surjective map such that $\varphi_{|T=0}=f_0\otimes \frac{A}{I}$.  
Then the question is, whether there is a
surjective 
map $F:P[T]\sur J$ such that (i) $F_{|T=0}=f_0$ and (ii) $F_{|Z}=\varphi$.
{\rm Assume $2n\geq d+3$.}
\end{question}

The statement  of  Question   \ref{homoConj}, is 
a simple translation
of the theorem of Nori \cite[\S 3 Appendix]{M2}, 
 on smooth vector bundles $V$ over smooth
manifolds $M$, using a vector bundle to projective module dictionary.
%
The Question  \ref{homoConj}, as stated, would fail to have an affirmative answer, 
without the regularity  hypothesis   \cite[Example 6.4]{BS1}.
 Even 
 when $A$ is regular, without the condition $2n\geq d+3$, the question would not have an affirmation answer (see \cite[Example 3.15]{BS1}). 
However,
existing results (see \cite{M2, BS1, BK}) indicate that with suitable  hypotheses the regularity  and/or transversality hypotheses may be spared.
Up to date, the best affirmative result   (assumes $2n\geq d+3$) on   (\ref{homoConj}) is due to Bhatwadekar and Keshari \cite{BK},
 preceded by \cite{MS, M2, MV, BS1}.

While the  Homotopy Question (\ref{homoConj})
always  had the flavor of being  central to the Homotopy Program,  it was never articulated as such. In fact, it was never well 
understood by the researchers how or why so?
 This article clarifies and establishes the centrality of the Homotopy Question (\ref{homoConj}).
%
The other half of these two pillars in this program is the definitions of Euler class groups. Followed by the outline given by Nori,
for integers $0\leq n \leq d$ and line bundles $L$, definitions of Euler class groups $E^n(A, L)$ were given in \cite{BS3, BS2, MY}. 
In fact, Nori originally outlined a definition of $E^d(A, A)$, when $A$ is regular (see \cite{MS}). For any projective $A$-module $P$, with $rank(P)=d$,
an Euler class $e(P)\in E^d(A, \wedge^dP)$ was defined and it was proved \cite{BS3} that 
$$
e(P)=0 \Llra P\cong Q\oplus A.
$$
When $rank(P) \leq d-1$,
a desire to define  a similar obstruction class $e(P)$,  in some appropriate obstruction group or set seemed too ambitious. 
We accomplish this goal, under additional conditions
(see Corollary \ref{e0e1Split}), by understanding the implicit Homotopy relations   in the statement of the Homotopy Question \ref{homoConj}. 
We introduce the following notations: 
$$
\left\{
\begin{array}{l}
{\mathcal LO}(P)= \left\{(I, \omega): \omega: P \sur \frac{I}{I^2},~{\rm is~a~surjective~map,~where}~I~{\rm is~an~ideal} \right\}\\
{\mathcal LO}^n(P)= \left\{(I, \omega)\in {\mathcal LO}(P): height(I)=n \right\}\\
{\mathcal LO}^n_c(P)= \left\{(I, \omega)\in {\mathcal LO}(P): height(I)=n,~{\rm and}~V(I)~{\rm is~connected} \right\}\\
\end{array}
\right.
$$
%
 There is a (chain) homotopy relation ingrained in the statement of (\ref{homoConj}), by substituting $T=0, 1$, on the set ${\mathcal LO}(P)$. 
The set of equivalence classes would be denoted by $\pi_0\left({\mathcal LO}(P)\right)$. 
In ${\mathcal LO}(P)$,
 there are two distinguished elements  $(0, 0), (A, 0) \in{\mathcal LO}(P)$,  
and their images in $\pi_0\left({\mathcal LO}(P)\right)$ are denoted, respectively, by ${\bf e}_0$ and ${\bf e}_1$. 
Define the obstruction class
\begin{equation}\label{eulerPDef}
\varepsilon(P):={\bf e}_0\in \pi_0\left({\mathcal LO}(P)\right),
\end{equation}
to be called the {\bf (Nori) Homotopy class} of $P$.
We give a summary of the    main results in this article, before making further introductory remarks.
Let $A$ and $P$ be as above. 
\bE
\item\label{spltCorIntro} (See Corollary \ref{e0e1Split}.) Suppose $A$ is essentially smooth, over an infinite perfect field $k$. Assume  $2n\geq d+3$, with $1/2\in k$. 
Then, we prove
$$
P\cong Q\oplus A   \Llra   \varepsilon(P)={\bf e}_1 ~~{\rm (" the~ additive ~zero";~see ~(\ref{monoiDIntro}))}.
$$
  
 \item\label{liftOrinetIntro} (See Theorem \ref{formalizeBSNori}.) Suppose $A$ is essentially smooth, over an infinite perfect field $k$, with $1/2\in k$.
  Assume  $2n\geq d+3$. Let $(I, \omega)\in {\mathcal LO}^n(P)$ and let 
 $[(I, \omega)]\in \pi_0({\mathcal LO}(P))$  be its image. Then,
 $\omega:P\sur \frac{I}{I^2}$  lifts 
 $$
{\rm  to ~a ~surjective~ map}\quad
 \Omega:P\sur I\quad \Llra\quad  \varepsilon(P)=[(I, \omega)]\in \pi_0({\mathcal LO}(P)).
 $$

\item\label{monoiDIntro}  (See Theorem \ref{abelianGroup}.) 
Assume $A$ is a regular ring, containing a field $k$, with $1/2\in k$. Assume $2n\geq d+2$. Then, we prove that $\pi_0({\mathcal LO}(P))$
has a natural structure of an abelian monoid. In this additive structure, ${\bf e}_1\in \pi_0({\mathcal LO}(P))$ is the identity. 
For $(I, \omega_1), (J, \omega_2)\in {\mathcal LO}^n(P)$, if $I+J=A$, the sum in $\pi_0({\mathcal LO}(P))$ is given by
$$
[(I, \omega_1)]+[(J, \omega_2)]= [(IJ, \omega_1\star \omega_2)]
$$
where $\omega_1\star \omega_2: P\sur \frac{IJ}{(IJ)^2}$ is obtained by combining $\omega_1$ and $\omega_2$, using Chinese remainder theorem.

Further, if $P=Q\oplus A$, then ${\bf e}_0={\bf e}_1$ and $\pi_0({\mathcal LO}(P))$ has a  structure of a group.

\item\label{eulerIntro}
To further establish centrality of  the Homotopy Question ( \ref{homoConj}) in this program, 
define Euler class group
$$
E(P):=\frac{{\BZ}({\mathcal LO}^n_c(P))}{{\SR}(P)}
$$
where ${\SR}(P)\subseteq {\BZ}({\mathcal LO}(P))$ is the subgroup generated by the  global orientations, namely, those $(I, \omega)\in {\mathcal LO}^n(P)$
such that, 
$\omega$ lifts to a surjective map $P\sur I$
({\it here $(I, \omega)$ is considered as an element in ${\BZ}({\mathcal LO}^n_c(P))$, by decomposing $I$ in to connected components}). 

\bE
\item\label{12SepSurMap} (See Definition \ref{EntopizeroNov14}.) Assume $A$ is a regular ring, containing a field $k$, with $1/2\in k$. 
Assume $2n\geq d+2$ and $P=Q\oplus A$. Then, we prove that there is a natural 
surjective group homomorphism
$$
\varphi: E(P) \sur \pi_0\left({\mathcal LO}(P) \right)
$$
\item\label{12SepIso} (See Theorem \ref{useBKNov14}.)  Further, assume $A$ is essentially smooth over an infinite perfect field $k$ and $1/2\in k$. If $2n\geq d+3$, 
then we prove that the homomorphism $\varphi$ is
an isomorphism.
\item (See Theorem \ref{oldGoldMethod}.) Assume $A$ is a noetherian commutative ring (without any regularity hypothesis), $P=Q\oplus A$ and 
$2n\geq d+3$. Let $(I, \omega)\in {\mathcal LO}^n(P)$. Assume its image 
$$
\overline{(I, \omega)}=0 \in E(P).
$$
Then, $\omega$ lifts to a surjective map $\Omega:P \sur I$.

\item(Corollary \ref{CorSallFour}.) In Section \ref{compChowSec}, we exploit the work on N.Mohan Kumar and M. M. Murthy \cite{MoM, Mu}, when the base field $k$ is algebraically closed, and 
$P\cong Q\oplus A$, with $rank(P)=n=d$. If $A$ is smooth and $1/2\in k$, we prove $\pi_0\left({\mathcal LO}(P)\right)\cong CH^d(A)$, where $CH^d(A)$ denotes the Chow group
of zero cycles.

\eE
\eE
The desire to define  an obstruction class, for $P$ 
to split off a free direct summand,  is age old and might have been  considered too bold. However, we are able to give such a definition (\ref{eulerPDef})  
 of an obstruction class $\varepsilon(P)$, and the result in item \ref{spltCorIntro} (Corollary \ref{e0e1Split}) establishes the 
splitting property. The result in item \ref{liftOrinetIntro} (Theorem \ref{formalizeBSNori}) was the main objective of the Homotopy Question (\ref{homoConj}), 
in such a homotopy obstruction theory set up. The structure of 
$\pi_0\left({\mathcal LO}(P)\right)$ has
 been an open problem since the inception of the Homotopy Program, while the exact nature  of the structure to expect was not clear. In item \ref{monoiDIntro}
(Theorem \ref{abelianGroup}) we settle this issue, by proving that the homotopy obstruction set 
$\pi_0\left({\mathcal LO}(P) \right)$ has structure of a monoid. The definition, in item \ref{eulerIntro}, of Euler class group 
$E(P)$  is new. Note, for a line bundle $L$, $E(L\oplus A^{n-1})$ coincides with $E^n(A, L)$, as defined in 
\cite{BS3, BS2, MY}. 
Further, under suitable smoothness conditions, the results in item \ref{eulerIntro}  (see \S \ref{secular} for more details),
 establish a relationship, as in (\ref{12SepSurMap}), (\ref{12SepIso}), between homotopy obstructions $\pi_0\left({\mathcal LO}(P)\right)$ and the Euler class group
$E(P)$, which ties together the two components of the germs of ideas originally given by Nori (around 1990). When $n=d$,  $k$ is algebraically closed, $P=Q\oplus A$,
under suitable other conditions, 
we establish that $\pi_0\left({\mathcal LO}(P)\right)$ coincides with the Chow group $CH^d(A)$ of zero cycles.

While we described our results above, in terms of $\pi_0\left({\mathcal LO}(P)\right)$, there are three other descriptions of $\pi_0\left({\mathcal LO}(P)\right)$ 
available in \S \ref{foundationSec}. 
We use these descriptions of $\pi_0\left({\mathcal LO}(P)\right)$ interchangeably.
Consider the notations:
$$
\left\{
\begin{array}{l}
\CQ(P)=\left\{(f, s)\in P^*\oplus A: s(1-s)\in f(P)  \right\}\\
\widetilde{\CQ}(P)=\{(f, p, s)\in P^*\oplus P\oplus A: f(p)+s(s-1)=0 \}\\
\widetilde{\CQ}'(P)=\{(f, p, z)\in P^*\oplus P\oplus A: f(p)+z^2=1 \}\\
\end{array}
\right.
$$
Given a polynomial extension $A \hra A[T]$, substituting $T=0, 1$, we have two set theoretic maps, in each case
$$
\left\{
\begin{array}{l}
\diagram 
\CQ(P) & \CQ(P[T]) \ar[r]^{T=1} \ar[l]_{T=0} & \CQ(P)
\enddiagram\\
\diagram 
\widetilde{\CQ}(P) & \widetilde{\CQ}(P[T]) \ar[r]^{T=1} \ar[l]_{T=0} & \widetilde{\CQ}(P)
\enddiagram\\
\diagram 
\widetilde{\CQ}'(P) & \widetilde{\CQ}'(P[T]) \ar[r]^{T=1} \ar[l]_{T=0} & \widetilde{\CQ}'(P)
\enddiagram\\
\end{array}
\right.
$$
These lead to  chain homotopy relations and accordingly, $\pi_0\left(\CQ(P) \right)$, $\pi_0\left(\widetilde{\CQ}(P)\right)$, $\pi_0\left(\widetilde{\CQ}'(P)\right)$ are defined. 
If and when $1/2\in A$ ({\it which we often assume}), there is a bijection 
$$
\widetilde{\CQ}(P)\iso \widetilde{\CQ}'(P),~{\rm which ~induces~a ~bijection}\quad 
\pi_0\left(\widetilde{\CQ}(P)\right)\iso\pi_0\left(\widetilde{\CQ}'(P)\right).
$$
In Section  \ref{foundationSec}, we  establish the following commutative diagram of natural bijections:
$$
\diagram
\pi_0\left(\widetilde{\CQ}(P)\right) \ar@{->>}[r]^{\overline{\nu}}_{\sim} \ar@{->>}[d]_{\overline{\eta}}^{\wr} 
& \pi_0\left(\CQ(P)\right)\ar@{->>}[dl]^{\overline{\eta}'}_{\sim}\\
\pi_0\left({\mathcal LO}(P)\right)&\\
\enddiagram
$$

We comment on the use of the phrase "Homotopy Program".
Perhaps, the phrase  was  first used by Mandal, in a conversation with Nori to describe this whole set of problems. 
Among what were encapsulated in the program are the following:
\bE
\item\label{Part1} ({\bf Part 1})
A coherent theory of obstructions, based on homotopy was expected. It was    also expected that these homotopy obstructions 
would come together with  the concept of Euler class groups. 
\item\label{PartTwo}({\bf Part 2}) The theory 
should reconcile with the ${\BA}^1$-homotopy  approach 
({\it also know as Motivic  or Chow-Witt group approach}
\cite{BM, Mo}). 
\item\label{Part3ree} ({\bf Part 3})
When $A$ is a real smooth  affine algebra, this algebraic homotopy obstruction theory should also reconcile with the   Topological counter part, in
the  sense analogous to \cite{MSh}.
\eE
Results in this article addresses Part 1 of this program, in a comprehensive manner.
In deed,  a coherent theory of homotopy obstructions is established, as was expected. 
Note that the theory is not expected to behave too well for the lower half of the range of $n=rank(P)$. 
When $P$ is not free, the definition of the Euler class group $E(P)$ is new, which was needed to bring the two components of the Homotopy 
Program together.
The destination of the road map that emerged out of the introduction of the Homotopy Question (\ref{homoConj})
was not very well understood. 
This article clarifies and brings us to that destination. This completes the Part 1 of the program.
This was accomplished entirely by the methods of  commutative algebra, which was possible due to the strength of the Homotopy Question.

In a sense, Part 2 of the program was resolved in  \cite{AF, MM1}
fairly  satisfactorily, by settling the problem of Fabien Morel \cite[pp. 13]{Mo} affirmatively.
However, 
this article
  reveals that the
Chow-Witt groups $\widetilde{CH}^n(A)\cong \pi_0\left(\widetilde{\CQ}(A^n)\right)$ 
correspond  only to the case $P=A^n$. 
This article raises newer questions (\ref{motivicQn}), what would be an appropriate ${\BA}^1$-homotopy interpretation for $\pi_0\left(\widetilde{\CQ}(P)\right)$, analogous to 
the original question of Morel \cite[pp.13]{Mo}.
 While such an ${\BA}^1$-homotopy interpretation would be of its own interest,
 this may become useful for 
further study of the structure of these monoids $\pi_0\left(\widetilde{\CQ}(P)\right)$, like finite generation and others.
The Part 3 of the program would have to be addressed  subsequently. It appears, there is no well formulated or well studied
topological counter part to these monoids $\pi_0\left({\mathcal LO}(P) \right)$, in the literature.

We contrast Euler class groups $E(A^n)=E^n(A, A)$, Chow Witt groups $\widetilde{CH}^n(A)$, and 
Homotopy obstructions $\pi_0\left({\mathcal LO}(P)\right)$, along with the history.
%
%
The goal, indeed a desire, that eventually emerged out of 
the introduction of the Homotopy Question
(\ref{homoConj}), along with the definition of $E^n(A,A)=E(A^{n})$, 
was to define an obstruction class $\varepsilon(P)$ in a suitable obstruction set (desirably a group), 
for $P$ to split off a free direct summand.  
However, the study of Euler class groups $E(A^{n})$ stole most of the attention, 
while the study of the implications of the Homotopy Question 
(\ref{homoConj}) was largely left ignored. 
The ${\BA}^1$-Homotopy approach to Euler class groups emerged out of the
 introduction of Chow-Witt groups $\widetilde{CH}^n(A)$  in 2000 \cite{BM}, followed by the book of Fabien Morel \cite{Mo}, in 2012. Morel defined a surjective map 
 $E^n(A, A) \sur  \widetilde{CH}^n(A)$, and indicated that this map could be an isomorphism \cite[pp. 13]{Mo}. Morel also established \cite[Chapter 8]{Mo}
 that $\widetilde{CH}^n(A)$ is in bijection with the so called naive Homotopy group, 
 that is $\pi_0\left({\mathcal LO}(A^n)\right)$. Using this description of $\widetilde{CH}^n(A)$, 
as envisioned by Morel, it was established in  \cite{AF, MM1}, that the map  $E^n(A^n) \sur  \widetilde{CH}^n(A)$ is an isomorphism, under suitable other hypotheses. 
Morel did not indicate that some variation of his ${\BA}^1$-Homotopy approach, may lead to a definition of an obstruction class $\varepsilon(P)$, to split off a free direct summand. 
Both $E(A^n)$ and $\widetilde{CH}^n(A)$ are invariants of $A$.
They are not precise enough to house such an obstruction class $\varepsilon(P)$.
Even the Euler class groups $E(P)$ defined in this article (\ref{eulerIntro})
 are not large enough to house such an obstruction $\varepsilon(P)$.
Nori provided the
precise insight (\ref{homoConj}), exactly what would work (around 1990). This article brings Nori's insight to the fullest fruition and establishes that 
the Homotopy obstruction set $\pi_0\left({\mathcal LO}(P)\right)$ is the appropriate set to house such an obstruction $\varepsilon(P)$.

We comment on the organization of this article. First and foremost, it is best that the reader is familiar with the introduction of \cite{MM1}.
In section \ref{foundationSec}, we lay out the basic definitions and the foundation of this article. In this section, we define the Homotopy 
obstruction set  $\pi_0\left({\mathcal LO}(P)\right)$, and give three other descriptions of 
 the same, as mentioned above. In section \ref{secHomotopy} we prove that the 
chain homotopy relations on  $\widetilde{{\CQ}}'(P)$, is indeed an 
equivalence relation, under further regularity hypotheses.
In section \ref{LiftingSec}, we prove  our main results on lifting and splitting, which are independent of  the additive structure on $\pi_0\left({\mathcal LO}(P) \right)$.
In section \ref{secInvolution}, 
we define the involution map $\Gamma: \pi_0\left(\widetilde{{\CQ}}(P)\right)\iso \pi_0\left(\widetilde{{\CQ}}(P)\right)$, 
which may be thought of as a substitute for the additive-inverse map, without any regard to the existence of any  additive structure 
on $\pi_0\left(\widetilde{{\CQ}}(P)\right)$.
In section \ref{secGroupStru}, we establish the monoid structure on $\pi_0\left(\widetilde{{\CQ}}(P)\right)$. 
In section \ref{secular}, we define the Euler class group $E(P)$, and compare it with the  homotopy obstruction monoid $\pi_0\left(\widetilde{{\CQ}}(P)\right)$,
as well with Chow group of zero cycles.
In the Appendix section \ref{motivicSec}, we define $\pi_0\left(\widetilde{{\CQ}}(P)\right): \Sch_A \lra  {\Sets}$, as pre-sheaf, and raise the question (\ref{motivicQn}) of its motivic interpretation.



\vspace{2mm}
\noindent{\bf Acknowledgement:} 
{\it The authors would like to thank Madhav V. Nori for his academic support, guidance and for sharing his invaluable insight on this program, over a long period of time. 
 }
\section{Foundation of  Homotopy Obstructions}\label{foundationSec}

In this section, we establish some notations and, for a  projective module $P$,
over a noetherian ring $A$, give several descriptions of the homotopy pre-sheaves.

\begin{notations}\label{nota}{\rm 
Throughout, $A$ will denote a commutative noetherian ring with $\dim A=d$ and $k$ will denote a field. Often, but not always, we will assume $1/2\in  A$
and/or $k\subseteq A$.

For  $A$-modules $M, N$, we denote $M[T]:=M\otimes A[T]$ and $M^*=Hom(M,A)$.
For
$f\in Hom(M,N)$, denote $f[T]:=f\otimes 1\in Hom(M[T], N[T])$. 
Homomorphisms  $f: M\lra \frac{I}{I^2}$ would be identified with the induced 
maps $\frac{M}{IM} \lra \frac{I}{I^2}$.

 For surjective homomorphisms $\omega_1: M\sur \frac{I_1}{I_1^2}$, $\omega_2: M\sur \frac{I_2}{I_2^2}$, where 
 $I_1, I_2$ are two ideals, with $I_1+I_2=A$,  $\omega_1\star\omega_2: M\sur \frac{I_1I_2}{(I_1I_2)^2}$ will denote the unique surjective map 
induced by $\omega_1, \omega_2$. 

For a projective $A$-module $P$,   $\BQ(P)=(\BQ(P), q)$ will denote the quadratic space  $\BH(P)\perp A$, where $\BH(P)=P^*\oplus P$ is the hyperbolic space.
So, $P^*\oplus P\oplus A$ is the underlying projective module of  $\BQ(P)$ and, for $(f, p, s)\in P^*\oplus P\oplus A$, $q(f, p, s)=f(p)+s^2$.

The category of (noetherian) schemes over $\spec{A}$ will be donated by $\Sch_A$.
The category of sets will be denoted by $\Sets$.
%
Given a pre-sheaf $\CF:\Sch_A \to \Sets$, and a scheme $X\in \Sch_A$, define $\pi_0(\CF)(X)$ by the pushout 
\begin{equation}\label{Defpi0CF}
\diagram 
\CF(X\times \BA^1) \ar[r]^{T=0}\ar[d]_{T=1} & \CF(X)\ar[d]\\
\CF(X) \ar[r] & \pi_0(\CF)(X)\\ 
\enddiagram 
\qquad {\rm in}~~\Sets
\end{equation}
So, $X \mapsto \pi_0(\CF)(X)$ is also a pre-sheaf on $ \Sch_A$.
 For  an affine scheme $X=\spec{B}\in \Sch_A$ and a pre-sheaf $\CF$, as above, we write $\CF(B):=\CF(\spec{B})$ and $\pi_0(\CF)(B):=\pi_0(\CF)(\spec{B})$.
 
 }
\end{notations}
Given a projective $A$-module $P$, 
 we define a homotopy obstruction set $\pi_0({\mathcal LO}(P))$ 
and establish various other descriptions of the same. These are analogous to similar obstruction sets available in the literature,
when $P=A^n$ is free.
%

%
\bD\label{defpizero2Oct}{\rm 
 Let $A$ be a noetherian commutative ring, $X=\spec{A}$ and $P$ be a projective $A$-module. 
 By a {\bf local $P$-orientation}, we mean 
a pair $(I, \omega)$ where $I$ is an ideal of $A$ and $\omega:P \sur \frac{I}{I^2}$ is a surjective homomorphism, which is 
identified with 
 surjective homomorphism $\frac{P}{IP} \sur \frac{I}{I^2}$, induced by $\omega$.  A local {\bf local $P$-orientation} will simply be 
 referred to as a {\bf local orientation}, when $P$ is understood.
%
Denote 
\begin{equation}\label{4LOetc}
\left\{
\begin{array}{l}
{\mathcal LO}(P)=\left\{(I, \omega): (I, \omega)~{\rm is~a~local}~P~{\rm orientation} \right\}\\
\CQ(P)=\left\{(f, s)\in P^*\oplus A: s(1-s)\in f(P)  \right\}\\
\widetilde{\CQ}(P)=\{(f, p, s)\in P^*\oplus P\oplus A: f(p)+s(s-1)=0 \}\\
\widetilde{\CQ}'(P)=\{(f, p, z)\in P^*\oplus P\oplus A: f(p)+z^2=1 \}\\
\end{array}
\right.
\end{equation}
There is a commutative diagram of set theoretic maps,  denoted as follows:
\begin{equation}\label{diaEtasNU}
\diagram
\widetilde{\CQ}(P) \ar@{->>}[r]^{\nu} \ar@{->>}[d]_{\eta} & \CQ(P)\ar@{->>}[dl]^{\eta'}\\
{\mathcal LO}(P)&\\
\enddiagram~~
{\rm where,~for}~(f, p, s)\in \widetilde{\CQ}(P), \quad \nu(f, p, s)=(f, s)\quad  
\end{equation}
and $\eta'(f, s)=\eta(f, p, s)= (I, \omega)$, where $I=f(P)+As$ and $\omega: P \sur \frac{I}{I^2}$
is the homomorphism
is induced by $f$. These maps $\eta, \eta', \nu$ are surjective. If and when $1/2\in A$ ({\it which we often assume}), there is also a bijection 
\begin{equation}\label{QPQprimeNov19}
\kappa: \widetilde{\CQ}(P)\iso \widetilde{\CQ}'(P)\quad{\rm sending}\quad 
(f, p, s) \mapsto \left(2f, 2p, 2s-1\right) 
\end{equation}
Now, suppose $P$ is a fixed projective $A$-module, and schemes $Y\in \underline{Sch}_A$ ,with $\pi: Y\to \spec{A}$.
Then, ${\mathcal LO}(\pi^*P)$, $\widetilde{\CQ}(\pi^*P)$, $\CQ(\pi^*P)$, $ \widetilde{\CQ}'(\pi^*P)$ are likewise defined (\ref{4LOetc}).
 The associations $Y \mapsto {\mathcal LO}(\pi^*P)$,
$Y \mapsto \CQ(\pi^*P)$,  $Y \mapsto \widetilde{\CQ}(\pi^*P)$, $Y \mapsto \widetilde{\CQ}'(\pi^*P)$ are pre-sheaves on $\Sch_A$. However, 
the pre sheaf nature of $Y \mapsto {\mathcal LO}(\pi^*P)$ requires some clarification. For example, for a ring homomorphism $\beta: A \lra B$, 
and $(I, \omega)\in {\mathcal LO}(P)$ is sent to $(\beta(P)B, \omega')\in {\mathcal LO}(P\otimes B)$, where $\omega': P\otimes B\lra \frac{\beta(I)B}{\beta(I^2)B}$
is induced by $\omega$.

By the pushout diagram (\ref{Defpi0CF}), applied to these pre-sheaves defines, the Homotopy obstructions pre-sheaves 
\begin{equation}\label{def4Hsheaf}
Y \mapsto \left\{
\begin{array}{l}
\pi_0\left({\mathcal LO}(P)\right)(Y),\\ 
\pi_0\left(\CQ(P)\right)(Y)\\
\pi_0\left(\widetilde{\CQ}(P)\right)(Y),\\ 
\pi_0\left(\widetilde{\CQ}'(P)\right)(Y). \\ 
\end{array}
\right.
\qquad {\rm are~defined.}
\end{equation}
%
 For historical reasons, we  explicitly define the  
Homotopy obstruction set
$\pi_0\left({\mathcal LO}(P)\right)$, by the 
pushout diagrams, in $\underline{Sets}$, as follows:
\begin{equation}\label{pizeroLOpush}
\diagram
{\mathcal LO}(P[T])\ar[r]^{T=0}\ar[d]_{T=1} &{\mathcal LO}(P)\ar[d]\\
{\mathcal LO}(P)\ar[r] & \pi_0\left({\mathcal LO}(P)\right)\\ 
\enddiagram 
\qquad  {\rm in}\quad \underline{Sets}.
\end{equation}
In deed, $\pi_0\left({\mathcal LO}(P)\right)$ was the Homotopy obstruction explicitly envisioned by Nori (see \cite{M2}). 
}
\eD


\vspace{2mm}
For the convenience of our discussions, we make the following notational adjustment.
\begin{notations}\label{notaAdjust}
{\rm
Until Section \ref{motivicSec}, we would  only be interested in  the value of 
the above homotopy pre sheaves (\ref{def4Hsheaf}),   when $Y=\spec{A}$. 
To simplify notations, we will make the following notational adjustment:
$$
\left\{
\begin{array}{l}
 \pi_0\left(\widetilde{{\CQ}}(P)\right):= \pi_0\left(\widetilde{{\CQ}}(P)\right)(A)\\
  \pi_0\left(\widetilde{{\CQ}}'(P)\right):= \pi_0\left(\widetilde{{\CQ}}'(P)\right)(A)\\
 \pi_0\left({\mathcal LO}(P)\right):= \pi_0\left({\mathcal LO}(P)\right)(A)\\
\end{array}
\right.
$$
Technically, as well, this adjustment would not make any difference. We would prove subsequently that all these sets are isomorphic, when $1/2\in A$.
This set $ \pi_0\left({\mathcal LO}(P)\right)$  would  be referred to as {\bf  Homotopy obstruction Set} of $P$. 
}
\end{notations}

We record, the following basic lemma 

\bL\label{basicMaps}
Use the notations as above  (\ref{notaAdjust})
and assume $1/2\in A$.
Then, the bijection $\kappa$, induces an isomorphism
$$
\overline{\kappa}: \pi_0\left(\widetilde{\CQ}(P)\right)\iso \pi_0\left(\widetilde{\CQ}'(P)\right)
$$
Further, the maps $\eta, \nu, \eta'$ (in diagram \ref{diaEtasNU})
 induce set theoretic maps, as denoted in the commutative diagram of  maps of pre-sheaves:
$$
\diagram
\pi_0\left(\widetilde{\CQ}(P)\right) \ar@{->>}[r]^{\overline{\nu}} \ar@{->>}[d]_{\overline{\eta}} 
& \pi_0\left(\CQ(P)\right)\ar@{->>}[dl]^{\overline{\eta}'}\\
\pi_0\left({\mathcal LO}(P)\right)&\\
\enddiagram
$$
\eL
\pf It follows from definition of pushout. 
$\eop$

\vspace{2mm}
We proceed to prove that, 
  the above is a commutative triangle of bijections: 
\begin{equation}\label{3obsWihBee}
\diagram
\pi_0\left(\widetilde{\CQ}(P)\right) \ar@{->>}[r]^{\overline{\nu}}_{\sim} \ar@{->>}[d]_{\overline{\eta}}^{\wr} 
& \pi_0\left(\CQ(P)\right)\ar@{->>}[dl]^{\overline{\eta}'}_{\sim}\\
\pi_0\left({\mathcal LO}(P)\right)&\\
\enddiagram
\end{equation}
We fix notations, for $(f, p, s)\in \widetilde{\CQ}(P)$, its equivalence class in $\pi_0\left(\widetilde{\CQ}(P)\right)$ will be denoted by 
$[(f, p, s)]$ and similar notations will be used for $(f, s)\in \CQ(P)$ and $(I, \omega)\in {\mathcal LO}(P)$.
Note, given $(I, \omega)\in \CL\CO(P)$, $\omega$ lifts to a homomorphism $f$, as follows:
\begin{equation}\label{aLiftOmega}
\diagram 
P \ar[r]^f\ar@{->>}[dr]_{\omega} & I \ar@{->>}[d]\\
& \frac{I}{I^2}\\
\enddiagram
\end{equation}
By Nakayama's lemma there is an element $s\in I$ such that $(1-s)I\subseteq f(P)$.
Consequently, $(f,s)\in \CQ(P)$ and $I=(f(P), s)$. This association would not be unique.
Such a pair $(f,s)\in \CQ(P)$ will be referred to as a lift of $(I,\omega)$ in $\CQ(P)$.
Now define an the 
 map:
\begin{equation}
\label{defOfZetaOct2}
\chi: \CL\CO(P) \lra \pi_0\left(\CQ(P) \right) \quad{\rm by~~~} \chi\left(I, \omega\right)=[(f, s)]\in \pi_0\left(\CQ(P) \right)
\end{equation}
where $(f, s) \in \CQ(P)$ is any lift of $(I,\omega)$ in $\CQ(P)$,  ({\it as in diagram \ref{aLiftOmega}}) and 
$[(f, s)]$ is its equivalence class.
In several lemmas, we establish that $\chi$ is well defined. 

\bL
\label{LindOfs2Oct} {\rm
Use the  notations as in (\ref{notaAdjust}). 
Let $(I, \omega_I)\in \CL\CO(P)$  and  $(f,s)\in \CQ(P)$ be a lift, as in diagram (\ref{aLiftOmega}). 
Further, assume that $t(1-t) \in f(P)$, with $I=(f(P), s)=(f(P), t)$.
%
Then 
$$
[(f, s)]=[(f, t)]\in  \pi_0\left(\CQ(P)\right).
$$
}
\eL
\pf 
First note, $(1-s)I \subseteq f(P)$ and $(1-t)I\subseteq f(P)$.
Write $I[T]=IA[T]$. So,
$$
I[T]=f(P)A[T]+sA[T]=f(P)A[T]+tA[T].
$$
Let   $S(T)=t+T(s-t)$. Clearly, $S(T)\in I[T]$. Further,
$$
{\bf Claim}: \qquad (1-S(T))I[T] \subseteq f(P)A[T] 
$$
We have $(1-S(T))I[T]=(1-S(T))[f(P)A[T]+sA[T]$. 
So, we only need to prove that $(1-S(T))s \in f(P)A[T]$. But
$$
(1-S(T))s= (1-t)s-T(s-t)s= (1-t)s+T[(s-1)s+(1-t)s] \in f(P)A[T] 
$$
So, the claim is established. Therefore, $(1-S(T))S(T)\in  f(P)A[T]$.
Denote $f[T]:=f\otimes 1_{A[T]}$. Then, $f[T]: P[T]\sur f(P)A[T]$ is a surjection.
Clearly, $(f[T], S(T))\in \CQ(P[T])$. Now,
$(f[T], S(T))_{T=0}=(f, t)$ and $(f[T], S(T))_{T=1}=(f, s)$. 
 \pic $\eop$  

\bL
\label{2Lifts2Oct}
Use the  notations as in (\ref{notaAdjust}). 
Suppose $(I, \omega)\in \CL\CO(P)$ and and $f, g$ be two lifts of $\omega$ as follows:
$$
\diagram 
P\ar@{->>}[r]^f\ar@{->>}[dr]_{\omega} & f(P)\ar@{->>}[d]\\ 
  & \frac{I}{I^2}\\
\enddiagram \quad {\rm and} \quad 
\diagram 
P\ar@{->>}[r]^g\ar@{->>}[dr]_{\omega} & g(P)\ar@{->>}[d]\\ 
  & \frac{I}{I^2}\\
\enddiagram
$$
$$
\ni ~~I=(f(P), s)=(g(P), t) ~{\rm and} ~s(1-s) \in f(P), ~t(1-t) \in g(P).
$$
Then 
$$
[(f, s)]=[(g, t)]\in  \pi_0\left(\CQ(P) \right)
$$
\eL
\pf Note,  $(g-f)(P)\subseteq I^2$. Let $F=f[T]+T(g[T]-f[T]) \in P[T]^*$. 
It is obvious that
$$
I[T]=F(P[T])+I[T]^2
$$
For completeness, we give a proof.
$$
\forall ~x\in I, x=(1-s)x+sx =f(p)+ sx \quad {\rm where} \quad p\in P~sx\in I^2
$$
So,
$$
{\rm (modulo ~}I[T]^2{\rm )}
\quad x \equiv f(p) \equiv F[T](p). 
$$
So,
$$
\exists \quad S(T)\in I[X]~\ni (1-S(T))I[T]\subseteq  F[T](P[T]) 
$$
So, $(F[T], S(T))\in \CQ(P[T])$. Therefore,
$$
[(f, S(0))]=[(F(0), S(0))]=[(F(1), S(1))]=[(g, S(1))] 
$$
Now, the proof is complete by (\ref{LindOfs2Oct}). $\eop$
\bT
Use the  notations as in (\ref{notaAdjust}). 
Let $(I,\omega)\in \CL\CO(P)$. Then, $\chi(I, \omega)$ as defined in equation (\ref{defOfZetaOct2}), is well defined.
\eT
\pf Follows from Lemma \ref{2Lifts2Oct}. $\eop$

\vspace{2mm}
Now, we prove that $\overline{\nu}$ is a bijection, as follows.
\bT\label{2pi0SamNov26}
Use the  notations as in (\ref{notaAdjust}).
Then, the map 
$$
\overline{\nu}: \pi_0\left(\widetilde{\CQ}(P) \right)\sur \pi_0\left(\CQ(P)\right)
\qquad {\rm is~a~bijection}. 
$$
\eT 
\pf 
Define a map 
 $\Psi_0: \CQ(P) \to \pi_0\left(\widetilde{\CQ}(P)\right)$ as follows:
 Given $(f,s)\in \CQ(P)$,  $\exists~p\in P~\ni ~f(p)=s(1-s)$. Define 
 $$
 \Psi_0(f, s):=[(f, p, s)]\in  \pi_0\left(\widetilde{\CQ}(P)\right).
 $$
We show that this association is a well defined map. To show this, suppose there is another 
$q\in P$ such that $f(q)=s(1-s)$. Note $f(p-q)=0$. So, $f[T](p+T(q-p))= f(p)+Tf(q-p)= f(p)+0=s(1-s)$.
Therefore,
$$
H(T):=\left(f[T], p+T(q-p), s) \right)\in \widetilde{\CQ}(P[T])
$$
and, hence
$$
H(0)=(f, p, s)\sim H(1)=(f, q, s).
$$
This establishes that $\Psi_0$ is well defined. Now, we   show that $\Psi_0$ is homotopy invariant. 
To see this, suppose $H(T)=(F, S(T))\in {\CQ}(P[T])$. Then, $S(T)(1-S(T))=F(p(T))$, for some 
$p(T)\in P[T]$. 
Write $\widetilde{H}=(F, p(T), S(T))\in \widetilde{\CQ}(P[T])$. 
So, 
$$
\Psi_0(F(0), s(0))=[\widetilde{H}(0)]= [\widetilde{H}(1)]= \Psi_0(F(1), S(1)) 
$$
This establishes that $\Psi_0$ factors through a map
$$
\Psi:  \pi_0\left(\CQ(P)\right) \to \pi_0\left(\widetilde{\CQ}(P)\right).
$$
It is easy to check that $\overline{\nu}$ and $\Psi$ are inverse of each other. 
\pic $\eop$ 


\vspace{2mm}
\bL\label{3reePINot}
{\rm 
Use the  notations as in (\ref{notaAdjust}).
Then, the 
 map $\chi: {\mathcal LO}(P) \lra \pi_0\left(\CQ(P)\right)$ (see (\ref{defOfZetaOct2})) 
 induces a well defined 
map $ \overline{\chi}: \pi_0\left({\mathcal LO}(P)\right) \lra \pi_0\left(\CQ(P)\right)$, 
 %
which is the inverse of the map $\overline{\eta}':\pi_0\left(\CQ(P)\right) \lra  \pi_0\left({\mathcal LO}(P)\right)$. 

Consequently, 
all the maps $\overline{\eta}$, $\overline{\eta}'$, $\overline{\nu}$ in diagram \ref{3obsWihBee},
are bijections.

}\eL
\pf The latter statement follows from the first one. 
Given a homotopy $H(T)\in {\mathcal LO}(P[T])$, it lifts to a homotopy $\widetilde{H}(T)=(F(T), S(T))\in {\CQ}(P[T])$.
So, $\chi(H(0))=[(F(0), S(0))]= [(F(1), S(1))]= \chi(H(1))$. So, $\chi$ is homotopy invariant, hence $\overline{\chi}$ 
is well defined. It is easy to see that this induced map in the inverse of $\overline{\eta}'$. 
%
%
%
\pic 
$\eop$


\bC
Use the  notations as in (\ref{notaAdjust}). Recall the notation 
$$
Q_{2n}(A)=\left\{(x_1, \ldots, x_n; y_1, \ldots, y_n; z)\in A^{2n+1}:\sum x_iy_i=z(1-z) \right\}.
$$
With $P=A^r=\oplus Ae_i$ is free,  then $Q_{2n}(A)\cong \widetilde{{\CQ}}(P)$  is a bijection.
This bijection induces a bijection $\pi_0\left(Q_{2n}\right)(A))\cong \pi_0\left(\widetilde{{\CQ}}(P)\right)$.
\eC

Before we proceed, we introduce the following notions.
\begin{notations}\label{notZetas}{\rm
Suppose $A$ is a commutative noetherian ring, with $\dim A=d$ and $P$ is a projective $A$-module, with $rank(P)=n$.
Denote $\zeta= \overline{\nu}^{-1}\chi:{\mathcal LO}(P) \lra \pi_0\left(\widetilde{{\CQ}}(P)\right)$ and $\zeta_0:\widetilde{{\CQ}}(P) \lra 
 \pi_0\left(\widetilde{{\CQ}}(P)\right)$. So, we have a commutative diagram:
 $$
 \diagram 
 \widetilde{{\CQ}}(P) \ar[rd]^{\zeta_0}\ar[d]_{\eta} & \\
 {\mathcal LO}(P)  \ar[r]_{\zeta\qquad} &  \pi_0\left(\widetilde{{\CQ}}(P)\right)\\
 \enddiagram
 $$
}\end{notations}

\vspace{2mm}
\begin{remark}{\rm 
The equation $\sum_{i=1}^nX_iY_i+Z(Z-1)=0$ would be the main motivation behind the definition of $\widetilde{{\CQ}}(P)$. 
For a field $k$ and the ring
${\SB}(k)=\frac{k[X_1, \ldots, X_n; Y_1, \ldots Y_n; Z]}{(\sum_{i=1}^nX_iY_i+Z(Z-1))}$, Swan \cite{Sw2} computed the total Chow ring $CH({\SB}(k))$ of ${\SB}(k)$.
This ring ${\SB}(k)$ is sometimes referred to as the {\bf universal ring}, for complete intersections. 
Using the structure of the Chow ring $CH({\SB}(k))$, together with Riemann-Roch Theorem,  Mohan Kumar and Nori \cite{Mk} proved that the 
ideal
 $I=(X_1, \ldots, X_n, Z){\SB}(k) \subseteq 
{\SB}(k)$,
 is not image of a projective ${\SB}(k)$-module of rank $n$.

In the more recent past, the notation $Q_{2n}:=\spec{{\SB}(k)}$ has been somewhat standard in the literature 
of the motivic approach to Euler class theory.

In fact, sometimes it would
 be  convenient to work with $\widetilde{{\CQ}}(P)$ than ${\mathcal LO}(P)$. This is due to the fact that, when $1/2\in A$,  $\widetilde{{\CQ}}(P)\cong  \widetilde{{\CQ}}'(P)$,
which has a nice quadratic structure that we can exploit (see \S \ref{secHomotopy}).
}
\end{remark}
\section{Homotopy Equivalence}\label{secHomotopy}

In this section, we prove the following key homotopy theorem.

\bT\label{PopThm3p4M2} 
Let $A$ be a regular ring over a  field $k$, with $1/2\in k$. Let $P$ be a projective $A$-module, with $rank(P)=n\geq 2$, 
and $(\BQ(P), q)=\BH(P) \perp A$
{\rm (see \ref{nota})}. Recall $\widetilde{{\CQ}}'(P)\subseteq \BQ(P)=P^*\oplus P \oplus A$.
Suppose $H(T)\in \widetilde{{\CQ}}'(P[T])$. Then, there is an orthogonal matrix $\sigma(T)\in O(\BQ(P), q)$, such that 
$$
H(T)=\sigma(T)(H(0))\qquad {\rm and}\qquad \sigma(0)=1.
$$
%
\eT
\pf 
Let $H(T)=(f(T), p(T), s(T)  \in \widetilde{{\CQ}}'(P[T])$ be a homotopy, as above. So, $H(0) \in \widetilde{{\CQ}}'(P)$. Then,
$$
A[T]H(T) \cong A[T]H(0) \cong (A[T], q_0) \qquad {\rm are ~isometric},
$$
where $q_0$ is the trivial quadratic space of rank one.
The bilinear inner product  in $\BQ(P)$ will be denoted $\langle -, -\rangle$.
We have the following split exact sequences of quadratic spaces:
$$
\diagram
0\ar[r] & K \ar[r] & \BQ(P[T]) \ar[rr]^{\langle H(T), -\rangle} & & A[T] \ar[r] & 0\\ 
0\ar[r] & K_0 \ar[r] & \BQ(P) \ar[rr]_{\langle H(0), -\rangle} & & A \ar[r] & 0\\ 
\enddiagram
$$
Therefore, 
$K=\left(A[T]H(T)\right)^{\perp}$,  $K_0=\left(AH(0)\right)^{\perp}$
are orthogonal complements. 
Write $\overline{K}:=K\otimes \frac{A[T]}{(T)}$. 
Note, for $\wp\in \spec{A}$, $\BQ(P)_{\wp} \cong (A, q_{2n+1})$, where $q_{2n+1}=\sum_{i=1}^nX_iY_i+Z^2$.
So, $\overline{K}_{\wp}\cong (K_0)_{\wp}$ are isometric. It is standard (see \cite[Lemma 4.1]{MM1}),
 that $(K_0)_{\wp}=(A_{\wp}H(0))^{\perp}\cong (A, q_{2n})_{\wp}$
where 
$q_{2n}=\sum_{i=1}^nX_iY_i$. In other words, $\overline{K}$ is locally trivial. By the the Quadratic version \cite[Theorem 3.5]{MM1} of  Lindel's 
theorem \cite{L},
there is an isometry $\tau:K \iso \overline{K}\otimes A[T]$. 
Further, it follows 
$\overline{K}=\left(RH(0)\right)^{\perp}\cong K_0$. 
Therefore, there is an isometry $\sigma_0:  \overline{K} \iso K_0$, which extends to an isometry
$\sigma_0\otimes 1: \overline{K}\otimes A[T] \iso K_0\otimes A[T]$.
Then,  $\sigma_1:=(\sigma_0\otimes 1)\tau: K \iso  K_0\otimes A[T]$ is an isometry.
Finally, note 
$$
(A[T]H(T), q_{|A[T]H(T)}) \cong (A[T], q_0) \cong (A[T]H(0), q_{|A[T]H(0)}).
$$
Now, consider the  diagram
\begin{equation}\label{linearMaps}
\diagram
0\ar[r] & K\ar[r]\ar[d]_{\sigma_1} & \BQ(P[T])  \ar@{-->}[d]^{\sigma(T)} \ar[rr]^{\langle H(T), -\rangle} & & A[T] \ar@{=}[d]\ar[r] & 0\\ 
0\ar[r] & K_0\otimes A[T] \ar[r] & \BQ(P[T])  \ar[rr]_{\langle H(0), -\rangle} & & A[T] \ar[r] & 0\\ 
\enddiagram
\end{equation}
of quadratic spaces. In this diagram, the horizontal lines are split exact sequences of quadratic spaces. Hence, there is an 
isometry $\sigma(T)\in O\left(A[T], q\right)$, such that the diagram commutes. That means, for all ${\bf v}\in  \BQ(P[T])$, we have
$\langle H(T), {\bf v}\rangle= \langle H(0), \sigma(T){\bf v}\rangle$. Replacing $\sigma(T)$ by $\sigma(T)^{-1}$, we have
 $\sigma(T)H(0)=H(T)$. So, we have $\sigma(0)H(0)=H(0)$. Again,
 by replacing $\sigma(T)$ by $\sigma(T)\sigma(0)^{-1}$, we have $\sigma(0)=1$.
\pic 
 $\eop$ 


\vspace{2mm}
The following Corollary would be of some importance for our future discussions.
\bC\label{directHomo}
Let $A$ be a regular ring over a  field $k$, with $1/2\in k$. Let $P$ be a projective $A$-module, with $rank(P)=n\geq 2$, 
and $(\BQ(P), q)=\BH(P) \perp A$.
Let ${\bf u}, {\bf v}\in \widetilde{{\CQ}}'(P)$ such that $[{\bf u}]= [{\bf v}]\in \pi_0\left(\widetilde{{\CQ}}'(P)\right)$.
Then, there is a homotopy $H(T)\in \widetilde{{\CQ}}'(P[T])$ such that $H(0)={\bf u}$ and $H(1)={\bf v}$.
Equivalently, for 
${\bf u}, {\bf v}\in \widetilde{{\CQ}}(P)$ if $\zeta_0({\bf u})= \zeta_0({\bf v})\in  \pi_0\left(\widetilde{{\CQ}}(P)\right)$, then
there is a homotopy $H(T)\in \widetilde{{\CQ}}(P[T])$ such that $H(0)={\bf u}$ and $H(1)={\bf v}$.
\eC
\pf Suppose ${\bf u}, {\bf v}\in \widetilde{{\CQ}}(P)$ such that $[{\bf u}]= [{\bf v}]\in \pi_0\left(\widetilde{{\CQ}}'(P)\right)$.
Then, there is a sequence of homotopies $H_1(T), \ldots , H_m(T) \in \widetilde{{\CQ}}'(P[T])$
such that ${\bf u}=:{\bf u}_0:=H_1(0)$, ${\bf u}_m:=H_m(1)={\bf v}$ and $\forall~i=1, \ldots, m-1$, we have ${\bf u}_i:=H_i(1)=H_{i+1}(0)$.
By Theorem \ref{PopThm3p4M2}, for $i=1, \ldots, m$ there are orthogonal matrices $\sigma_i(T)\in O(\BH(P[T]), q)$ such that
$\sigma_i(0)=1$ and 
$H_i(T)=\sigma_i(T)H_i(0) =\sigma_i(T){\bf u}_{i-1}$. Therefore, ${\bf u}_i=H_i(1) = \sigma_i(1){\bf u}_{i-1}$.

Write $H(T)=\sigma_m(T)\cdots \sigma_1(T){\bf u}_0$. Then, $H(T)\in \widetilde{{\CQ}}'(P[T])$ and $H(0)= {\bf u}_0$ and $H(1)= {\bf u}_m$.
This establishes first part of the statement on $\pi_0\left(\widetilde{{\CQ}}'(P)\right)$. The latter assertion on $\pi_0\left(\widetilde{{\CQ}}(P)\right)$ follows from the former, 
by the bijective correspondences 
$\widetilde{{\CQ}}'(P)  \iso \widetilde{{\CQ}}(P)$ and $\widetilde{{\CQ}}'(P[T]) \iso \widetilde{{\CQ}}(P[T])$. 
\tcp $\eop$
\begin{remark}\label{equivHomRelDec2}{\rm
Another way to state (\ref{directHomo}) would be  that the homotopy relation on $\widetilde{{\CQ}}(P)$ is actually an equivalence relation.
}
\end{remark}

In a slightly more formal language, the above is summarized as follows.
\bT\label{summHomT}
Let $A$ be a regular ring over a  field $k$, with $1/2\in k$. Let $P$ be a projective $A$-module, with $rank(P)=n\geq 2$, 
and $(\BQ(P), q)=\BH(P) \perp A$.
For, $\sigma(T)\in O\left(\BQ(P), q\right)$ and ${\bf u}\in \widetilde{{\CQ}}'(P)$, define the  {\rm (left)} action 
$\sigma(T){\bf u}:=\sigma(1){\bf u}\in \widetilde{{\CQ}}'(P)$. 
Denote $O\left(\BQ(P), q, T\right)=\{\sigma(T)\in O\left(\BQ(P), q\right): \sigma(0)=1 \}$.
Then, the map
$$
\frac{\widetilde{{\CQ}}'(P)}{O\left(\BQ(P), q, T\right)} 
\lra \pi_0\left(\widetilde{{\CQ}}'(P)\right) \quad {\rm is~a~bijection}.
$$
\eT
\pf Similar to the proof of (\ref{directHomo}). $\eop$
\section{Homotopy Triviality and Lifting} \label{LiftingSec}
In this section, under further smoothness  conditions, we establish that for $(I, \omega_I)\in {\mathcal LO}(P)$, 
 the triviality of  $\zeta(I, \omega_I)$ implies that $\omega_I$ lifts to a surjective map $P\sur I$.
We start this section with the following notations and definitions.
\bD\label{basePoint}{\rm
Suppose $A$ is a commutative noetherian ring, with $\dim A=d$ and $P$ is a projective $A$-module, with $rank(P)=n$.
There are two distinguished points in $\widetilde{{\CQ}}(P)$, namely:
$$
{\bf 0}:=(0, 0, 0) \in  \widetilde{{\CQ}}(P), \quad {\bf 1}:=(0, 0, 1) \in  \widetilde{{\CQ}}(P)
$$
$$
{\rm We~ denote}\quad
{\bf e}_0= \zeta_0({\bf 0}) \in \pi_0\left(\widetilde{{\CQ}}(P)\right), \quad {\rm and}\quad 
{\bf e}_1= \zeta_0({\bf 1}) \in \pi_0\left(\widetilde{{\CQ}}(P)\right).
$$
Use the same notations ${\bf e}_0, ~{\bf e}_1\in \pi_0\left({\mathcal LO}(P)\right)\cong \pi_0\left(\widetilde{{\CQ}}(P)\right)$, 
to denote their respective  images. 
Define the obstruction class
$$
\varepsilon(P):= {\bf e}_0 \in \pi_0\left({\mathcal LO}(P)\right)\cong \pi_0\left(\widetilde{{\CQ}}(P)\right).
$$
In the light of (\ref{homoConj}), $\varepsilon(P)$ will be referred to as {\bf (Nori) Homotopy  Class} of $P$, 
which may sometimes be shortened.
Note, for any $f\in P^*$ and $p\in P$,
$\varepsilon(P):={\bf e}_0= \zeta_0(f, 0, 0)= \zeta_0(0, p, 0) \in \pi_0\left(\widetilde{{\CQ}}(P)\right)$.

}
\eD
%

We record the following obvious observation.
\bL\label{Bibekananda}
Suppose $A$ is a commutative noetherian ring with $\dim A=d$ and $P$ is a projective $A$-module. 
Let $p\in P$ and $f\in P^*$ be such that $f(p)=1$ {\rm (i. e. $P\cong Q\oplus A$)}. Let 
$$
{\bf 0}=(0,0,0), ~
{\bf u}=(f, 0, 0), ~{\bf 1}=(0, 0, 1) \in \widetilde{{\CQ}}(P).
$$
Then, $\zeta_0({\bf 0})=\zeta_0({\bf u})
= \zeta_0({\bf 1}) \in \pi_0\left(\widetilde{{\CQ}}(P)\right)$. In other words,
$$
\varepsilon(P)={\bf e}_0={\bf e}_1.
$$
\eL
\pf The first equality is obvious and was mentioned above (\ref{basePoint}). To prove the second equality,
write $H(T)=\left((1-T)f, Tp, T \right)$. Then, $(1-T)f(Tp)=T(1-T)$. So, $H(T)\in \widetilde{{\CQ}}(P[T])$.
We have $H(0)= {\bf u}$ and $H(1)= (0, p, 1)$. 

Now write $G(T)=(0, (1-T)p, 1) )\in \widetilde{{\CQ}}(P[T])$.  Then, $G(0)= (0, p, 1)$ and $G(1)= (0, 0, 1)$.
\pic $\eop$
%


\vspace{2mm}
The following is the main result in this section. 
%
\bT\label{formalizeBSNori}
Suppose $A$ is an essentially smooth ring over an infinite perfect field $k$, with $1/2\in k$ and $\dim A=d$. 
Let $P$ be a projective $A$-module with $rank(P)=n$, with $2n\geq d+3$. 
Let ${\bf 0}=(0, 0, 0)\in \widetilde{{\CQ}}(P)$, as in (\ref{basePoint}).
%
Suppose $(I, \omega_I)\in {\mathcal LO}(P)$, with $height(I)\geq n$. Then, $\omega_I$ lifts to a surjective map $P\sur I$ if and only if 
$\varepsilon(P)=\zeta(I, \omega_I)$. 
\eT
\pf 
Suppose $\omega_I$ lifts to a surjective map $f:P \sur I$.   
Write $H(T)= (f(T), 0, 0)\in \widetilde{{\CQ}}(P[T])$. Then, $\zeta(I, \omega_I)=\zeta_0(H(1))=\zeta_0(H(0))=  \zeta_0({\bf 0})=\varepsilon(P)$.

Conversely, suppose $\zeta(I, \omega_I)=\zeta_0({\bf 0})$.
For notational convenience, fix $f_0\in P^*$,   and let ${\bf v}_0=(f_0, 0, 0)\in \widetilde{{\CQ}}(P)$. Then, $\zeta(I, \omega_I)=\zeta_0({\bf 0})=\zeta_0({\bf v}_0)$.  
There is an element ${\bf u}=(f_1, p_1, s_1) \in \widetilde{{\CQ}}(P)$ such that $\eta({\bf u})=(I, \omega_I)$.
By Moving Lemma argument \ref{movingLemm} (below), we can assume that  $height( f_0(P))\geq n$ and $height( f_1(P))\geq n$.
We have, $\zeta_0({\bf u})=\zeta_0({\bf v}_0)$. By (\ref{directHomo}), there is a homotopy $H(T)=(f(T), p(T), S(T))\in \widetilde{{\CQ}}(P[T])$ such that 
$H(0) = {\bf v}_0$ and $H(1)={\bf u}$. Write $\eta(H(T))=(J, \Omega)$.  We would apply \cite[Theorem 4.13]{BK}, for which we would need $height(J)
\geq n$. So, we modify $H(T)$, as follows. Denote $Z(T)=1-S(T)$. Write 
$\SP=\{\wp\in \spec{A[T]}: height(\wp)\leq n-1, T(1-T)Z(T)\notin \wp \}$.
 Define a generalized dimension $\delta:\SP \lra \BN$ by $\delta(\wp)=height(\wp)$ for $\wp\in \SP$.
 Then,
$\forall~\wp\in \SP$, we have $\delta(\wp)\leq n-1$.
Now, $(f(T), T(1-T)Z(T)^2)\in P[T]^*\oplus A[T]$ is basic on $\SP$. 
 So, there is an element $g(T)\in P[T]^*$ such that
 $F(T)=f(T)+T(1-T)Z(T)^2g(T)$
  is basic on $\SP$.
 It follows, $F(0)=f(0)$ and $F(1)=f(1)$.
 
We have $Z(T)(1-Z(T))=(1-s(T))s(T)=$
$$
f(T)(p(T))=F(T)(p(T))-T(1-T)Z(T)^2g(T)(p(T))
$$
Write $\CJ=(f(T)(P[T]), Z(T))$. Then $\CJ=(F(T)(P[T]), Z(T))$. Write $M=\frac{\CJ}{F(T)(P[T])}$. Let $p_1, \ldots, p_m$ be a set of 
generators of $P$. So, $\CJ$ is generated by the images of $f(T)(p_1), \ldots, f(T)(p_m), Z(T)$. 
Use "overline" to denoted images in $M$ and repeat the proof of Nakayama's Lemma, as follows:
$$
\left(
\begin{array}{c}
\overline{f(T)(p_1)}\\
\overline{f(T)(p_2)}\\
\cdots\\
\overline{f(T)(p_m)}\\
\overline{Z(T)}\\
\end{array}
 \right)
 =
 \left(
\begin{array}{ccccc}
0 & 0 & \cdots & 0 & -T(1-T)Z(T)g(T)(p_1)\\
0 & 0 & \cdots & 0 & -T(1-T)Z(T)g(T)(p_2)\\
0 & 0 & \cdots & 0 & \cdots\\
0 & 0 & \cdots & 0 & -T(1-T)Z(T)g(T)(p_m)\\
0 & 0 & \cdots & 0 & Z(T)-T(1-T)Z(T)g(T)(p(T))\\
\end{array}
 \right)
 \left(
\begin{array}{c}
\overline{f(T)(p_1)}\\
\overline{f(T)(p_2)}\\
\cdots\\
\overline{f(T)(p_m)}\\
\overline{Z(T)}\\
\end{array}
 \right)
$$
So,
$$
 \left(
\begin{array}{ccccc}
1 & 0 & \cdots & 0 & T(1-T)Z(T)g(T)(p_1)\\
0 & 1 & \cdots & 0 & T(1-T)Z(T)g(T)(p_2)\\
0 & 0 & \cdots & 0 & \cdots\\
0 & 0 & \cdots & 1 & T(1-T)Z(T)g(T)(p_m)\\
0 & 0 & \cdots & 0 &1-Z(T)+ T(1-T)Z(T)g(T)(p(T))\\
\end{array}
 \right)
 \left(
\begin{array}{c}
\overline{f(T)(p_1)}\\
\overline{f(T)(p_2)}\\
\cdots\\
\overline{f(T)(p_m)}\\
\overline{Z(T)}\\
\end{array}
 \right)
 =
 \left(
\begin{array}{c}
0\\
0\\
\cdots\\
0\\
0\\
\end{array}
 \right)
 $$
With $Z'(T)=Z(T)-T(1-T)Z(T)g(T)(p(T))$, the determinant of this matrix is $1-Z'(T)$. It follows, $(1-Z'(T))\CJ\subseteq F(T)(P[T])$.
So,\\
$(1-Z'(T))Z'(T)=F(T)(q(T))$ for some $q(T)\in P[T]$. Note, $Z'(0)=Z(0)$ and $Z'(1)=Z(1)$.
Therefore, $(F(T), q(T), Z'(T))\in \widetilde{{\CQ}}(P[T])$. Also, with
  $S'(T)=1-Z'(T)$, $(F(T), q(T), S'(T))\in \widetilde{{\CQ}}(P[T])$. We have  \\
$S'(T)(1-S'(T))=(1-Z'(T))Z'(T)=F(T)(q(T))$\\
$S'(0)=1-Z'(0)=1-Z(0)=S(0)=0$ and\\
$S'(1))=1-Z'(1)=1-Z(1)=S(1)$.\\
Write $\CH(T)= (F(T), q(T), S'(T))$ and $\eta(\CH(T))=(J', \Omega')$. It is clear $\CH(0)=(f_0, q(0), 0)$, $\CH(1)=(f_1, q(1), S(1))$.
So, $\eta(\CH(0))=\eta({\bf v}_0)$ and $\eta(\CH(1))=\eta({\bf u})=(I, \omega_I)$.

 We have $J'=(F(T)(P[T]), S'(T))$. We claim that $height(J')\geq n$. To see this,
 let $J'\subseteq \wp\in \spec{A[T]}$. If $T\in \wp$, then $I_0 \subseteq \wp$ and hence $height(\wp)\geq n$. Likewise, 
  if $1-T\in \wp$, then $I \subseteq \wp$ and hence
 $height(\wp)\geq n$. So, we assume $T(1-T)\notin \wp$. 
 If $Z(T)\in \wp$, then $\CJ=(F(T)(P[T]), Z'(T))= (F(T)(P[T]), Z(T))\subseteq \wp$, which is impossible because $S'(T)\in \wp$.
 So, $T(1-T)Z(T)\notin \wp$. Since $F$ is basic on $\SP$, $height(\wp)\geq n$. This establishes the claim.
 
 So, $\CH(T)=(F(T), q(T), S'(T))\in \widetilde{{\CQ}}(P[T])$ is such that $\eta(\CH(0))=(I_0, \omega_{I_0})$, 
 $\eta(\CH(1))=(I, \omega_{I})$ and with $\eta(\CH(T))=(J', \Omega')$, we have $height(J') \geq n$. 
 If $T\in \wp \in Ass\left(\frac{A[T]}{J'}\right)$ then $(J'(0), T)=(I_0, T)\subseteq \wp$. Then, $height(\wp)\geq n+1$.
 This is impossible because $A[T]$ is regular (Cohen-Macaulay) and $J'$ is local complete intersection ideal.
 Hence, 
 $$
 \diagram
 \frac{A[T]}{TJ'}\ar[r]^{T=0}\ar[d] & A\ar[d]\\
 \frac{A[T]}{J'}\ar[r] & \frac{A}{J'(0)}\\
 \enddiagram
 $$
 is a patching diagram (see (\ref{patchBK}) below). So, the map $\Omega': P[T] \sur \frac{J'}{(J')^2}$ and $f_0:P\sur I_0$ combines to give a surjective 
 maps $\phi:P[T] \sur \frac{J'}{T(J')^2}$. Now, by \cite[Theorem 4.13]{BK}, there is a surjective homomorphism 
 $\varphi:P[T] \sur J'$ such that $\varphi(0)=f_0$ and $\varphi\otimes \frac{A[T]}{J'}=\Omega'$. Now, it follows that 
 $\varphi(1)$ is a lift of $\omega_I$.
 \tcp $\eop$





\vspace{2mm}
We used the following lemma above, while it needs a proof. The standard references for Patching diagrams are \cite{Mi, R, O2}. 
We will be specific in the following statement, because
the 
literature does not seem complete regarding definitions Patching diagrams of modules that are not projective. 
\bL\label{patchBK}
Let $R$ be a noetherian commutative ring and
 $A=R[T]$. Take $J=AT$ and  $I$ is locally complete intersection ideal of height $r$ and $T:\frac{A}{I} \hra \frac{A}{I}$ is injective 
(i. e.
$T\notin \wp \in Ass\left( \frac{A}{I}\right)$). Then,
$$
\diagram 
\frac{I}{TI^2} \ar[r]\ar[d] & \frac{I}{TI}\ar[d]\\
\frac{I}{I^2} \ar[r] & \frac{I}{I^2+TI}\\
\enddiagram
$$
is a Patching diagram, in the sense that it is a Cartesian square.
Further,
\bE
\item $\frac{I}{TI} \iso I(0)$.
\item $\frac{I}{I^2+TI}\iso \frac{I(0)}{I(0)^2}$.
\eE
\eL
\pf The patching diagram follows, 
 because 
 $I^2\cap (TI)=TI^2$.
 
 To see this, first
we have $TI^2 \subseteq I^2\cap (TI)$. Suppose $f\in I^2\cap (TI)$. Then, $f=Tg$ with $g\in I$. Now, consider the map
$$
T: \frac{I}{I^2} \lra \frac{I}{I^2}
$$
Since $\frac{I}{I^2}$ projective and $T:\frac{A}{I} \hra \frac{A}{I}$ is injective, $T$ is also injective on $\frac{I}{I^2}$. So, $g\in I^2$.
So, $f=Tg\in TI^2$. 

Now, we prove $\frac{I}{TI} \iso I(0)$. Obviously, the map is on to. Suppose $f(T)\in I$ and $f(0)=0$.Then, $f=Tg$. Since $T$ is non zero divisor on 
$\frac{A}{I}$, $g\in I$. So, $f\in TI$.

Finally, we prove $\frac{I}{I^2+TI}\iso \frac{I(0)}{I(0)^2}$. Again, the map is on to. Suppose $f(T)\in I$ and $f(0)\in I(0)^2$.
Then, $f(0)=\sum f_i(0)g_i(0)$. Then, $f-\sum f_ig_i\in (T)\cap I=TI$ (by the above, if we like). So, $f\in I^2+TI$.

$\eop$

We close this section with the following  "moving lemma argument", which  is fairly standard. A number of  variations of the same (\ref{movingLemm}) 
would be among the frequently used tools for the rest of our discussions.
\bL[Moving Lemma]\label{movingLemm}
Suppose $A$ is a commutative noetherian ring with $\dim A=d$ and $P$ is a projective $A$-module with $rank(P)=n$.
Assume $2n\geq d+1$.
Let $K\subseteq A$ be an ideal with $height(K)\geq n$ and $(I, \omega_I)\in {\mathcal LO}(P)$. 
Then, there is an element ${\bf v}=(f, p, s)\in \widetilde{{\CQ}}(P))$ such that 
$\eta({\bf v}) =(I, \omega_I)$. Further, with $J=f(P)+A(1-s)$, we have $height(J)\geq n$ and $J+K=A$.
\eL
\pf 
Let $f_0:P \sur I$ be any lift of $\omega_I$. Then, $I=f_0(P)+I^2$. By Nakayama's Lemma, 
there is an element $t\in I$,
such that $(1-t)I \subseteq f_0(P)$. Therefore, $t(1-t)=f_0(p_0)$ for some $p_0\in P$.
({\it Readers are referred to \cite{M1} regarding generalities on Basic Element Theory and generalized
dimension functions}.)
Write 
$$
\SP=\left\{\wp\in \spec{A}: t\notin \wp, ~{\rm and~either}~K\subseteq \wp~{\rm or}~height(\wp)\leq n-1 \right\}
$$
There is a generalized dimension function (see \cite{M1})  $\delta: \SP \lra \BN$, such that $\delta(\wp)\leq n-1~\forall~\wp\in \SP$.
Now $(f_0, t^2)\in P^* \oplus A$ is basic on $\SP$. So, there is an element $g\in P^*$ such that $f:=f_0+t^2g$ is basic on $\SP$.
It follows, $f(P)+At=f_0(P)+At=I$ and $I=f(P)+I^2$. 
By Nakayama's Lemma, 
there is an element $s\in I$,
such that $(1-s)I \subseteq f(P)$ and hence $f(p)=s(1-s)$, for some $s\in I$. Hence, $I=(f(P), s)$.
Now, write $J=f(P)+A(1-s)$. For $J \subseteq \wp\in \spec{A}$, $s\notin \wp$ and hence $t\notin \wp$.
Since, $f$ is basic on $\SP$, $height(\wp)\geq n$. This establishes, $height(J)\geq n$. 

Now suppose  $J+K \subseteq \wp\in \spec{A}$. By the same argument above, $t\notin \wp$. Hence,
$\wp\in \SP$. This is Impossible, because $f$ is basic on $\SP$. So,  $J+K=A$. Now, ${\bf v}=(f, p, s)
\in \widetilde{{\CQ}}(P)$, satisfies the requirement.
$\eop$



\vspace{2mm}
The following is a converse of Lemma \ref{Bibekananda}.
\bC\label{e0e1Split}
Suppose $A$ is an essentially smooth ring over an infinite perfect  field $k$, with $1/2\in k$ and $\dim A=d$. 
Let $P$ be a projective $A$-module with $rank(P)=n$. Assume $2n\geq d+3$. Then, 
$$
\varepsilon(P)={\bf e}_1 \quad \Llra \quad P\cong Q\oplus A
$$
for some projective $A$-module $Q$.
\eC
\pf Suppose $P\cong Q\oplus A$. Then, by (\ref{Bibekananda}), 
$\varepsilon(P)={\bf e}_0 ={\bf e}_1$. Conversely, suppose $\varepsilon(P)={\bf e}_0 ={\bf e}_1$. Fix 
$f_0\in P^*$ such that $height(f_0(P))=n$. Then, $\zeta_0(f_0, 0, 0)={\bf e}_0 ={\bf e}_1$. Then, it follow from Theorem 
\ref{formalizeBSNori} that $\eta(0, 0, 1)$ lifts to a surjective map $P\sur A$. \tcp $\eop$

\section{The Involution}\label{secInvolution}
 
In this section,
we introduce an involution map  $\Gamma: \pi_0\left(\widetilde{{\CQ}}(P)\right) \lra \pi_0\left(\widetilde{{\CQ}}(P)\right)$. 
This can be thought of as a substitute to additive  inverse map, without any regard to existence of an addition.

\bD\label{invoDefi}
Suppose $A$ is a commutative ring and $P$ is a projective $A$-module, with $rank(P)=n$.
For $(f, p, s) \in \widetilde{{\CQ}}(P)$, define $\Gamma(f, p, s)= (f, p, 1-s)$.
This 
 association, ${\bf v}\mapsto \Gamma({\bf v})$, establishes a bijective correspondence
 $$
 \Gamma:  \widetilde{{\CQ}}(P)\iso  \widetilde{{\CQ}}(P), \quad{\rm such ~that}\quad \Gamma^2=1.
 $$
 We would say that $\Gamma$ is an involution  
 on $\widetilde{{\CQ}}(P)$, which will be a key instrument in the subsequent discussions. {\rm (This notation $\Gamma$ will be among the standard notations throughout this article.)}
\eD
We record the following  obvious lemma.
\bL\label{obvious}
Suppose $A$ is a commutative ring and $P$ is a projective $A$-module, with $rank(P)=n$ and 
 $\Gamma: \widetilde{{\CQ}}(P)\iso  \widetilde{{\CQ}}(P)$ is the involution.
Let 
 ${\bf v}=(f, p, s)\in \widetilde{{\CQ}}(P)$ and denote $\eta({\bf v})=(I, \omega_I)$ and $\eta(\Gamma({\bf v}))=(J, \omega_J)$.
Then,
\bE
\item $I\cap J=f(P)$.
 
\item For $H(T)\in \widetilde{{\CQ}}(P[T])$, we have
 $\Gamma(H(T))_{T=t}= \Gamma(H(t))$.
\item\label{obNov1216} Therefore,
$
\forall~{\bf v}, {\bf w}\in \widetilde{{\CQ}}(P)
\qquad  
\zeta_0({\bf v}) =\zeta_0({\bf w}) \Llra  \zeta_0(\Gamma({\bf v})) =\zeta_0(\Gamma({\bf w})).
$
\eE
\eL


\vspace{2mm}
In deed, $\Gamma$ factors through an involution on $ \pi_0\left(\widetilde{{\CQ}}(P)\right)$, as follows.
 
\bC\label{theInvolution} 
Suppose $A$ is a commutative ring and $P$ is a projective $A$-module, with $rank(P)=n$.
Then, the involution $\Gamma:  \widetilde{{\CQ}}(P)\iso  \widetilde{{\CQ}}(P)$ induces a bijective map
$\widetilde{\Gamma}:  \pi_0\left(\widetilde{{\CQ}}(P)\right) \iso  \pi_0\left(\widetilde{{\CQ}}(P)\right)$, such that 
$\widetilde{\Gamma}^2=1$ and $\zeta_0\Gamma= \widetilde{\Gamma}\zeta_0$. We say $\widetilde{\Gamma}$  is  an involution.
{\rm (The notation $\widetilde{\Gamma}$ will also be among our standard notations throughout this article.)}
\eC 
\pf  First, consider the map $\zeta_0\Gamma: \widetilde{{\CQ}}(P) \lra \pi_0\left(\widetilde{{\CQ}}(P)\right)$.
For, $H(T)\in \widetilde{{\CQ}}(P[T])$, we have $\zeta_0\Gamma(H(0))= \zeta_0\Gamma(H(1))$. Therefore, $\zeta_0\Gamma$
is homotopy invariant. Hence, it induces  a well defined map 
$\widetilde{\Gamma}:  \pi_0\left(\widetilde{{\CQ}}(P)\right) \iso  \pi_0\left(\widetilde{{\CQ}}(P)\right)$. 
Clearly, $\widetilde{\Gamma}^2=1$ and $\widetilde{\Gamma}$ is a bijection.
\pic $\eop$

\vspace{2mm}
The following is a way to compute the involution.
\bC\label{altF_i4involution}
Suppose $A$ is a commutative ring and $P$ is a projective $A$-module, with $rank(P)=n$.
Suppose $(I, \omega)\in {\mathcal LO}(P)$. For any ${\bf v}=(f, p, s) \in \widetilde{{\CQ}}(P)$ with 
$\eta({\bf v})= (I, \omega)$, write $\eta(\Gamma({\bf v})) =(J, \omega_J)$.
Then, 
$$
\widetilde{\Gamma}(\zeta(I, \omega))= \zeta(J, \omega_J)\in \pi_0\left(\widetilde{{\CQ}}(P)\right).
$$ 
\eC
 \pf Obvious.
$\eop$
 

The following is another version of the Moving Lemma \ref{movingLemm}. 
\bL[Moving Representation]\label{movingRep}
Suppose $A$ is a commutative ring, with $\dim A=d$.
Let  $P$ be a projective $A$-module, with $rank(P)=n$ and $2n\geq d+1$.
Let $x\in \pi_0\left(\widetilde{{\CQ}}(P)\right)$ and let $K\subseteq A$ be an ideal with $height(K)\geq n$. Then, there 
is a local $P$-orientation $(J, \omega_J)\in {\mathcal LO}(P)$ such that $x=\zeta(J, \omega_J)$, $height(J)\geq n$ and $J+K=A$.
\eL
\pf 
Let $x=\zeta(I, \omega_I)$. First, $\eta({\bf u})=(I, \omega_I)$ for some ${\bf u}\in \widetilde{{\CQ}}(P)$.
 Denote   $(I_0, \omega_{I_0}):=\eta(\Gamma({\bf u}))$. Then, 
$\tilde{\Gamma}(x)=\zeta(I_0, \omega_{I_0})$.

Now, we apply Moving Lemma \ref{movingLemm}, to $(I_0, \omega_{I_0})$ and $K$. 
There is ${\bf v}\in \widetilde{{\CQ}}(P)$, such that $\eta({\bf v})=(I_0, \omega_{I_0})$, and 
with $\eta(\Gamma({\bf v}))=(J, \omega_J)$, we have $height(J)\geq n$ and $J+K=A$. Now, $x=\tilde{\Gamma} (\tilde{\Gamma} (x))
= \tilde{\Gamma} (\zeta(I_0, \omega_{I_0}))=\zeta(J, \omega_J)$.
\pic $\eop$

\section{The Monoid Structure on $\pi_0\left({\mathcal LO}(P)\right)$} \label{secGroupStru}

In this section, we  define and establish a natural monoid
 structure on the homotopy obstruction set $\pi_0\left({\mathcal LO}(P)\right)\cong \pi_0\left(\widetilde{{\CQ}}(P)\right)$ ,
  when $2rank(P)\geq \dim A+2$ and $A$ is a regular ring 
over a  field $k$, with $1/2\in k$. We start with the following basic 
ingredient of the group structure. 
\bD\label{psudoPlus}
Let $A$ be a commutative noetherian ring, with $\dim A=d$,
 and $P$ be a projective $A$-module, with $rank(P)=n\geq 2$.
 Let
$(I, \omega_I), (J, \omega_J)\in {\mathcal LO}(P)$ be such that $I+J=A$. 
Let $\omega:=\omega_I\star \omega_J:  P \sur \frac{IJ}{(IJ)^2}$ be the unique surjective map induced by $\omega_I, \omega_J$.
We define a pseudo-sum 
$$
(I, \omega_I) \hat{+} (J, \omega_J) := (IJ, \omega)\in \pi_0\left({\mathcal LO}(P)\right).
$$
Note, pseudo-sum commutes.
\eD 
In the rest of this section, we establish that the pseudo sum respects homotopy, when $2n\geq d+2$, and 
$A$ is a regular ring 
over a  field $k$, with $1/2\in k$. Consequently, this leads to a addition operation on $\pi_0\left({\mathcal LO}(P)\right)$.
The following is the key lemma.

.



\bL\label{HLK1hoK2}
Let $A$ be a commutative noetherian ring
 and $P$ be a projective $A$-module, with $\dim A=d$, $rank(P)=n$, and $2n\geq d+2$.
Consider a homotopy
 $$
 H(T)=(f(T), p(T), Z(T))\in \widetilde{{\CQ}}(P[T]).
 $$
  Write $\eta(H(0))=
(K_0, \omega_{K_0})$ and $\eta(H(1))=(K_1, \omega_{K_1})$.
Further suppose $(J, \omega_J)\in {\mathcal LO}(P)$ such that $K_0+J=K_1+J=A$ and $height(J)\geq n$. Then, there is a homotopy 
$\CH(T)\in \widetilde{{\CQ}}(P[T])$ such that $\eta(\CH(0))=(K_0J, \omega_{K_0J})$ and $\eta(\CH(1))=(K_1J, \omega_{K_1J})$,
where, for $i=0, 1$ $\omega_{K_iJ}:=\omega_{K_i}\star \omega_J: P \sur \frac{K_iJ}{(K_iJ)^2}$.
Consequently,
$$
(K_1, \omega_{K_1}) \hat{+} (J, \omega_J) =(K_2, \omega_{K_2}) \hat{+} (J, \omega_J)\in \pi_0\left({\mathcal LO}(P)\right).
$$

\eL
\pf  We will write $f=f(T)$, $p=p(T)$ and $Z=Z(T)$. Dnote $Y=1-Z$ and $\eta(\Gamma(H(T))= (\BJ, \omega_{\BJ})$. Then, $\BJ= (f(P[T]), Y)$.
Write 
$$
\SP=\{\wp\in \spec{A[T]}: YT(1-T)\notin \wp,  J\subseteq \wp \}.
$$
 There is a generalized dimension function $\delta: \SP \lra \BN$ such that 
 $\forall~\wp\in \SP$,
 $\delta(\wp)\leq \dim\left(\frac{A[T]}{JA[T]}\right) \leq d+1-height(J) \leq d+1-n\leq n-1$.
 Further, $(f, Y^2T(1-T))$ is a basic element in $P[T]^*\oplus A[T]$, on $\SP$. Therefore, there is an element
$\lambda:=\lambda(T) \in P[T]^*$ such that
$$
f'=f+Y^2T(1-T)\lambda
~~ {\rm is~basic~on}~\SP.\quad 
{\rm So,} \quad  f'(0)=f(0), ~f'(1)=f(1).
$$
 We have $\BJ=(f(P[T]), Y)= (f'(P[T]), Y)$. Further,
$$
Z(1-Z)=Y(1-Y) =f(p)= f'(p) - Y^2T(1-T)\lambda(p).
$$
So,
$$
Y= f'(p) - Y^2T(1-T)\lambda(p) +Y^2
$$
 Write $M=\frac{\BJ}{f'(P[T])}$.  Let $p_1, \ldots, p_m$ be a set of generators of $P$.
 Use "overline" to indicate images in 
$M$. We intend to repeat the proof of Nakayama's Lemma and we have 
$$
\left(\begin{array}{c}
\overline{f(p_1)}\\ \overline{f(p_2)}\\ \cdots \\ \overline{f(p_m)}\\   \overline{Y}\\
\end{array} \right)
=
\left(\begin{array}{ccccc}
0 & 0 & \cdots & 0 & -\lambda(p_1)YT(1-T)\\
0 & 0 & \cdots & 0 & -\lambda(p_2)YT(1-T)\\
\cdots &\cdots &\cdots &\cdots &\cdots\\
0 & 0 & \cdots & 0 & -\lambda(p_m)YT(1-T)\\
0 & 0 & 0 & 0 & Y-  \lambda(p)YT(1-T) \\
\end{array} \right)
\left(\begin{array}{c}
\overline{f(p_1)}\\ \overline{f(p_2)}\\ \cdots \\ \overline{f(p_m)}\\   \overline{Y}\\
\end{array} \right)
\Lra 
$$
$$
\left(\begin{array}{ccccc}
1 & 0 & \cdots & 0 & \lambda(p_1)YT(1-T)\\
0 & 1 & \cdots & 0 & \lambda(p_2)YT(1-T)\\
\cdots &\cdots &\cdots &\cdots &\cdots\\
0 & 0 & \cdots & 1 & \lambda(p_m)YT(1-T)\\
0 & 0 & 0 & 0 & 1-Y+ \lambda(p)YT(1-T) \\
\end{array} \right)
\left(\begin{array}{c}
\overline{f(p_1)}\\ \overline{f(p_2)}\\ \cdots \\ \overline{f(p_m)}\\   \overline{Y}\\
\end{array} \right)
=
\left(\begin{array}{c}
0 \\ 0 \\ \cdots \\ 0\\ 0 \\
\end{array} \right)
$$
Multiplying by the adjoint matrix and computing the determinant, with $Y'= Y- \lambda(p)YT(1-T)$, we have
$$
(1-Y')\BJ\subseteq f'(P[T]).
$$
We have 
$Y'(0)=Y(0)=1-Z(0)$, $Y'(1)=Y(1)=1-Z(1)$.
Further, 
$$
Y'(1-Y')=f'(p') 
 \qquad{\rm for~some~polynomials}\quad p'\in P[T].
$$
$$
{\rm Therefore}\qquad H'(T)= (f'. p', Y) \in \widetilde{{\CQ}}(P[T]).
$$
We have
$$
\BJ= (f(P[T], Y)= (f'(P[T]), Y)=(f'(P[T]), Y').
$$
In fact, $\eta(H'(T))= (\BJ, \omega_{\BJ})$ and write $\eta(\Gamma(H'(T)))= (\BI, \omega_{\BI})$.   
Claim
$$
\BI+JA[T] =A[T].\qquad i.e. \qquad (f'(P[T]), 1-Y')+JA[T] =A[T].
$$
To see this, let 
$$
\BI+JA[T]\subseteq \wp \in \spec{A[T]}
$$
\bE
\item If $Y\in \wp$ then $\BJ=(f'(P[T]), Y) =(f'(P[T]), Y') \subseteq \wp$.
So, $Y'\in \wp$, which is 
impossible, since $1-Y'\in \wp$. 
So, $\wp \in D(Y)$.
\item 
Since $f'$
is unimodular of $\SP$ and $\wp \in D(Y)$ , we must have $T(1-T)\in \wp$.
\item Now, $T\in \wp$ implies, 
$$
\BI(0)+J=(f'(0)(P), 1-Y'(0))+J= (f(0)(P), 1-Y(0))+J= K_0+J=A\subseteq \wp
$$
%
%
which is impossible. 
\item Likewise, $1-T\in \wp$ implies, 
$$
\BI(1)+J= 
(f'(1)(P), 1-Y'(1))+J= (f(0)(P), 1-Y(1))+J= K_1+J=A\subseteq \wp.
$$
This is also impossible.
\eE
This establishes the claim.
%
Recall, $\omega_{\BI}: P[T] \sur \frac{\BI}{\BI^2}$ is induced by $f'$. 
Extend  $\omega_J:A^n \sur \frac{J}{J^2}$ to a surjective map  $\omega_{JA[T]}:A[T]^n \sur \frac{JA[T]}{J^2A[T]}$.
Let 
$$
\Omega:=\omega_{\BI}\star\omega_{JA[T]}  : P[T]\sur \frac{J\BI}{J^2\BI^2} \quad {\rm be~induced~by}~\Omega_{\BI}, ~{\rm and}~\omega_{JA[T]}.
$$
Now, there is a homotopy$\CH(T)\in \widetilde{{\CQ}}(P[T])$, such that  $\eta(\CH(T))= (\BI JA[T], \Omega)$. Specializing at $T=0$ and $T=1$, we have
$$
\eta(\CH(0))= (K_0J, \omega_{K_0J}), \qquad \eta(\CH(1))= (K_1J, \omega_{K_1J}).
$$
The proof is complete. $\eop$


\vspace{2mm}
Now, we define addition on $\pi_0\left({\mathcal LO}(P) \right)$.
\bD\label{directlyAdd}
Let $A$ be a regular, containing a field $k$, with $1/2\in k$, with $\dim A=d$,
Let $P$ be a projective $A$-module, with $rank(P)=n\geq 2$, and $2n\geq d+2$.
Let $x, y\in \pi_0\left({\mathcal LO}(P)\right)$. By Moving Lemma \ref{movingLemm}, we can write $x=[(I, \omega_I)]$, $y=[(J, \omega_J)]$, for some 
$(I, \omega_I), (J, \omega_J)\in{\mathcal LO}(P)$,
 with $height(IJ)\geq n$, 
 and $I+J=A$. Define 
$$
x+y:=(I, \omega_I)\hat{+} (J, \omega_J) \in \pi_0\left({\mathcal LO}(P)\right)
\qquad {\rm as~defined~in~ (\ref{psudoPlus})}.
$$
{\rm We establish that $x+y$ is well defined (\ref{wellShortDef}).}
\eD

\vspace{2mm} 
\bP\label{wellShortDef}
Under the setup and notations, as in (\ref{directlyAdd}), $x+y$ is well defined.
\eP
\pf Let $x=[(I_1, \omega_{I_1})], y=[(J_1, \omega_{J_1})]\in \pi_0\left({\mathcal LO}(P)\right)$, be another pair of choices, as in 
(\ref{directlyAdd}). That means, 
 $height(I_1J_1)\geq n$,  $I_1+J_1=A$. We   prove
$$
 (I, \omega_I) \hat{+} (J, \omega_J)=  (I_1, \omega_{I_1}) \hat{+} (J_1, \omega_{J_1}).
 $$
By Moving Lemma \ref{movingLemm}, there is $(K, \omega_K) \in {\mathcal LO}(P)$ such that 
$x=[(K, \omega_K)]$, $hieight(K)\geq n$ and  $K+ I_1\cap J_1=A$.

 We have ${\bf u}, {\bf u}_1\in\widetilde{\CQ}(P)$ such that $\eta({\bf u})= (I, \omega_I)$, and $\eta({\bf u}_1)= (K, \omega_{K})$.
Since $x=[(I, \omega_I)]=[(K, \omega_{K})] \in \pi_0\left({\mathcal LO}(P)\right)$, it follows ${\bf u}$, ${\bf u}_1$ are equivalent in 
$\widetilde{\CQ}(P)$. By (\ref{directHomo}), there is homotopy $H(T) \in \widetilde{{\CQ}}(P[T])$ such that $H(0) = {\bf u}$, and $H(1)={\bf u}_1$. 
It follows from
 Lemma \ref{HLK1hoK2}, 
 $$
 (I, \omega_I) \hat{+} (J, \omega_J)=  (K, \omega_{K}) \hat{+} (J, \omega_J)= (J, \omega_J) \hat{+} (K, \omega_{K})
 $$
 Likewise, the above 
 $$ 
 =
 (J_1, \omega_{J_1}) \hat{+} (K, \omega_{K})=  (K, \omega_{K})  \hat{+} (J_1, \omega_{J_1})
 =  (I_1, \omega_{I_1})  \hat{+} (J_1, \omega_{J_1})
 $$
\pic $\eop$
The   final statement on the binary structure on $\pi_0\left({\mathcal LO}(P)\right) \cong \pi_0\left(\widetilde{\CQ}(P)\right)$, 
is as follows.
\bT\label{abelianGroup}
Suppose $A$ is a regular ring over a field $k$, with $1/2\in k$ and $\dim A=d$. 
Let $P$ be a projective $A$-module with $rank(P)=n$. Assume $2n\geq d+2$.
%
{\rm (Subsequently, we use the notations in $\pi_0\left({\mathcal LO}(P)\right)$ and $\pi_0\left(\widetilde{\CQ}(P)\right)$ interchangeably.)}
Then, the addition operation on $\pi_0\left({\mathcal LO}(P)\right)$, 
defined in (\ref{directlyAdd})
has the following properties. 
\bE
\item\label{MonOne} The addition in $\pi_0\left({\mathcal LO}(P)\right)$ 
is commutative and associative.
Further, the image ${\bf e}_1:=[(A, 0)]\in \pi_0\left({\mathcal LO}(P)\right)$, of $(0, 0, 1)\in \widetilde{\CQ}(P)$,  acts as the additive identity in $\pi_0\left({\mathcal LO}(P)\right)$
In other words, $\pi_0\left({\mathcal LO}(P)\right)$ has a structure of an abelian monoid.
\item\label{GaisInv} 
Recall the involution map $\widetilde{\Gamma}: \pi_0\left({\mathcal LO}(P)\right) \iso \pi_0\left({\mathcal LO}(P)\right)$. 
Let  ${\bf e}_0:=[(0, 0)]\in \pi_0\left({\mathcal LO}(P)\right)$ be the image  $(0, 0, 0)\in \widetilde{\CQ}(P)$. Then,  $\forall ~x\in \pi_0\left({\mathcal LO}(P)\right)$ 
$x+\widetilde{\Gamma}(x)={\bf e}_0$.
\item \label{abGrpIFF} 
If ${\bf e}_0={\bf e}_1\in \pi_0\left({\mathcal LO}(P)\right)$, then
$\pi_0\left({\mathcal LO}(P)\right)$
 is an abelian group, under this addition. 
 {\rm (Recall (\ref{e0e1Split}), if $2n\geq d+3$, and if $A$ is essentially smooth over an infinite perfect field, then  ${\bf e}_0={\bf e}_1$ if and only if $P\cong Q\oplus A$.)}
\eE
\eT
\pf  
Given $x, y, z\in \pi_0\left({\mathcal LO}(P)\right)$, by   the Moving Lemma \ref{movingLemm}, we can write 
$$
x=[(K, \omega_K)],~y=[(I, \omega_I)]~z=[(J, \omega_J)]\quad \ni~K+I=K+J=I+J=A
$$
and $height(K)\geq n$, $height(I)\geq n$, $height(J)\geq n$. 
By definition (\ref{directlyAdd}),
$$
(x+y)+z=((K, \omega_K)\hat{+}(I, \omega_I))\hat{+} (J, \omega_J)=x+(y+z).
$$
$$
{\rm and} \qquad \qquad 
x+y= (K, \omega_K)\hat{+}(I, \omega_I)= (I, \omega_I)\hat{+}(K, \omega_K)=y+x.
$$
So, the associativity and commutativity hold. 
It is  obvious that, for all $x\in \pi_0\left({\mathcal LO}(P)\right)$, we have  $x+{\bf e}_1=x$. So, ${\bf e}_1$ acts as the additive identity.
This establishes (\ref{MonOne}).

Let $x=[(K, \omega_K)] \in \pi_0\left({\mathcal LO}(P)\right)$, with $height(K)\geq n$.
There is ${\bf u}=(f, p, s) \in \widetilde{{\CQ}}(P)$, with $\eta({\bf u})= (K, \omega_K)$.
Write 
$\eta(\Gamma({\bf u}))=(I_1, \omega_{I_1})$.
We can assume $height(I_1) \geq n$. 
 It follows. 
$$
x+\widetilde{\Gamma}(x)=\zeta_0(f, 0, 0)={\bf e}_0.\qquad {\rm This~ establishes ~(\ref{GaisInv}).}
$$

 If ${\bf e}_0={\bf e}_1$, it follows from  (\ref{GaisInv})
that, $\pi_0\left({\mathcal LO}(P)\right)$ has a group structure. This establishes  (\ref{abGrpIFF}).

 \tcp $\eop$


%


\begin{remark} \label{wittpiZeroDef}{\rm 
Use the notation as in (\ref{abelianGroup}). When ${\bf e}_0\neq {\bf e_1}$, the results in (\ref{abelianGroup}) describe a situation 
similar to  the construction of Witt group, from the monoid of isometry classes quadratic spaces. 

For  $x, y \in \pi_0\left({\mathcal LO}(P)\right)$ define $x\sim y$ if $x+n{\bf e}_0= y+m{\bf e}_0$, for integers $m, n\geq 0$.
This is easily checked to be an equivalence relation. 
Let $\BE\left(\pi_0\left({\mathcal LO}(P)\right)\right)$ be the set of all equivalence classes. Then, 
$\BE\left(\pi_0\left({\mathcal LO}(P)\right)\right)$ has a structure of an abelian group, induced by that additive structure 
on $\pi_0\left({\mathcal LO}(P)\right)$.
The natural  map
$$
\ell: \pi_0\left({\mathcal LO}(P)\right)
\sur \BE\left(\pi_0\left({\mathcal LO}(P)\right)\right)
$$
is a surjective homomorphism of monoids. 
The identity element of $\BE\left(\pi_0\left({\mathcal LO}(P)\right)\right)$ is  $\ell({\bf e}_0)=\ell({\bf e}_1)$.
For $x\in \pi_0\left({\mathcal LO}(P)\right)$, the additive inverse of $\ell(x)$ is $\ell(\widetilde{\Gamma}(x))$.

Clearly, if ${\bf e}_0= {\bf e_1}$, then $\BE\left(\pi_0\left({\mathcal LO}(P)\right)\right)=\pi_0\left({\mathcal LO}(P)\right)$.

} 

\end{remark}

 \section{The Euler Class Groups}\label{secular}

Suppose $A$ is a noetherian commutative ring with $\dim A=d$ and $P$ is a projective $A$-module, with $rank(P)=n$.
In this section, in analogy to the definition of the Euler class groups $E^n(A)$ in \cite{BS2, MY}, we define a group 
$E(P)$, which would also be called the Euler class group of $P$. Subsequently, we compare $E(P)$ with $\pi_0\left({\mathcal LO}(P)\right)$. 
%
%
 Also, refer to some superfluous aspect of the definitions 
in \cite{BS2, MY}, pointed out in \cite{MM1}. 
({\it In the sequel, for a set $S$, the free abelian group generated by $S$ will be denoted by} $\BZ(S)$).
\bD\label{newEuler}
{\rm
Suppose $A$ is a noetherian commutative ring, with $\dim A=d$ and $P$ is a projective $A$-module, with $rank(P)=n\geq 0$.
Denote,
$$
\left\{
\begin{array}{ll}
{\mathcal LO}^n(P)=\{(I, \omega_I)\in {\mathcal LO}(P):  height(I)= n \},\\
{\mathcal LO}_c^n(P)=\{(I, \omega_I)\in {\mathcal LO}(P):  V(I)~{\rm is~connected~and}~height(I)= n \}.\\
\end{array}
\right.
$$
Let $(I, \omega_I)\in {\mathcal LO}^n(P)$ and $I=\cap_{i=1}^mI_i$ be a decomposition, where $V(I_i)\subseteq 
\spec{A}$ are connected. 
The local orientation $(I, \omega_I) \in {\mathcal LO}^n(P)$ induce  $(I_i, \omega_{I_i})\in {\mathcal LO}^n_c(P)$, for $i=1, \ldots, m$.
Denote 
 $$
 (I, \omega_I)_{{\BZ}}=\sum_{i=1}^m(I_i, \omega_{I_i})\in \BZ\left({\mathcal LO}^n_c(P)\right).
 $$
A local orientation $(I, \omega_I)\in {\mathcal LO}(P)$ would be called {\bf global}, if $\omega_I$ lifts to a surjective map $P\sur I$.
 Let $\SR(P)$ denote the subgroup of $\BZ\left({\mathcal LO}^c(P)\right)$, generated by the set 
$\left\{(I, \omega_I)_{{\BZ}}:   (I, \omega_I)\in {\mathcal LO}^n(P),~{\rm and~is~global} \right\}$.

Define
$$
E(P)= \frac{\BZ\left({\mathcal LO}^n_c(P)\right)}{\SR(P)}
\quad {\rm to~be~called~the~Euler~class~group~of}~~P.
$$
Images of $(I, \omega_I)\in {\mathcal LO}(P)$ in $E(P)$ will be denoted by $\overline{(I, \omega_I)}$, which is same as the image of $(I, \omega_I)_E$.

}
\eD

Subsequently, we assume ${\bf e}_0={\bf e}_1 \in \pi_0\left({\mathcal LO}(P)\right)$, and hence $\pi_0\left({\mathcal LO}(P)\right)$ is a group.
In this case,   we define 
 a  homomorphism  $\rho: E(P) \lra \pi_0\left({\mathcal LO}(P)\right)$, 
 as follows.

\bD\label{EntopizeroNov14}{\rm 
Suppose $A$ is a regular  ring over a field $k$, with $1/2\in k$ and  $\dim A=d$, and $P$ is a projective $A$-module with $rank(P)=n$.
Assume  $2n\geq d+2$.   ({\it Use the notations in }(\ref{wittpiZeroDef})).
The  restriction $\beta$, of the map  $\zeta$, to ${\mathcal LO}^n_c(P)$,  gives the following commutative diagram:
$$
\diagram
{\mathcal LO}^n_c(P)\ar[dr]^{\beta} \ar@{^(->}[d]& \\
{\mathcal LO}^n(P)\ar@{->>}[r]_{\zeta} &\pi_0\left({\mathcal LO}(P)\right)\\ 
\enddiagram
$$
 We assume ${\bf e}_0={\bf e}_1$. So, $\pi_0\left({\mathcal LO}(P)\right)$ 
 has the structure of an abelian group. The  map $\beta$ extends to   
  group homomorphism 
    $\rho_0: \BZ\left({\mathcal LO}^n_c(P)\right)\lra\pi_0\left({\mathcal LO}(P)\right)$. 

  
  Now suppose $(I, \omega_I) \in {\mathcal LO}^n(P)$ be global. Let $f:P\sur I$ be a lift of $\omega$ and $I=\cap_{i=1}^m I_i$ be a decomposition of $I$ in to connected components. Then, 
  $$
  (I, \omega)_{{\BZ}}=\sum_{i=1}^m(I_i, \omega_i) \in \BZ\left({\mathcal LO}^n_c(P)\right).
  $$
  We have 
  $$
  \rho_0 \left((I, \omega)_{{\BZ}}\right)=
\sum_{i=1}^m[(I_i, \omega_i)]  = [\eta(f, 0, 0))]= {\bf e}_0 ={\bf e}_1
  $$
  Therefore, $\rho_0$ factors through a 
  group homomorphism $\rho: E(P)\sur \pi_0\left({\mathcal LO}(P)\right)$.
In fact, 
  $\rho$ is surjective. 
}
\eD 
\pf 
We only need to give a proof that $\rho$ is surjective.
For $x\in \pi_0\left({\mathcal LO}(P)\right)$, by   \ref{movingLemm}, $x=[(I, \omega_I)]$ for some 
$(I, \omega_I)\in {\mathcal LO}^n(P)$. Let $I =\cap_{i=1}^mI_i$ be a decomposition, with $V(I_i)$ connected and $\omega_i:
P\sur \frac{I_i}{I_i^2}$ be the surjective map induced by $\omega_I$. Then,
$$
\rho_0((I, \omega_I)_{{\BZ}})=\sum_{i=1}^m  [(I_i, \omega_i)]= [(I, \omega_I)]\in \pi_0\left({\mathcal LO}(P)\right).
$$
So, $\rho_0$ is surjective and hence so is $\rho$.
\tcp $\eop$

%

\bT\label{useBKNov14}
Suppose $k$ is an infinite perfect field, with $1/2\in k$ and $A$ is an essentially smooth ring over $k$, with $\dim A=d$. 
Suppose $P$ is a projective $A$-module with $rank(P)=n$
and  $2n\geq d+3$. Assume $P\cong Q\oplus A$.
Then, $\pi_0\left({\mathcal LO}(P)\right)$ is an abelian group and 
the homomorphism 
$\rho: E(P) \lra \pi_0\left({\mathcal LO}(P)\right)$ 
 is an isomorphism.
\eT
\pf We only need to prove that $\rho$ is injective. Let $\rho(x)=0$ for some $x\in E(P)$. We can write 
$x=\overline{(I, \omega_I)}$,
for some $(I, \omega_I) \in {\mathcal LO}^n(P)$. 
By Lemma \ref{Bibekananda}, we have $[(I, \omega_I)]= {\bf e}_1={\bf e}_0$. 
It follows from Theorem \ref{formalizeBSNori} that $\omega_I$ lifts to a surjective map $f:P\sur I$.
Therefore,
$(I, \omega_I)$  is global. Hence $x=\overline{(I, \omega_I)}=0$.
So, $\rho$ is an isomorphism.
 \tcp $\eop$
  
\bC\label{noCostEulerNov25}
Suppose $k$ is an infinite perfect field, with $1/2\in k$ and $A$ is an essentially smooth ring over $k$, with $\dim A=d$. 
Suppose
$P$ is a projective $A$-module with $rank(P)=n$
and  $2n\geq d+3$. Assume $P\cong Q\oplus A$. Suppose $(I, \omega_I)\in {\mathcal LO}^n(P)$ and 
$\overline{(I, \omega_I)}=0\in E(P)$. Then, $\omega_I$ lifts to a surjective homomorphism $P\sur I$.
\eC
 \pf It is immediate from Theorem \ref{useBKNov14}. $\eop$

\vspace{2mm}

In fact, a  stronger version (\ref{BS2Prop3p3Nov25}) of (\ref{noCostEulerNov25}) follows,    by   the same arguments as  in 
\cite{BS2}, .
%

 \subsection{The Vanishing of Euler cycles}
 We use the notations as in Definition \ref{newEuler}. An element $x\in E(P)$  is, sometimes,   referred to as an Euler cycle. In this subsection,
 we prove a less restrictive version of Corollary \ref{noCostEulerNov25}.
 We will follow the arguments in the proof of \cite[Theorem 4.2]{BS2}, 
 which mainly depends on the availability of Subtraction and Addition Principles.
   Accordingly, the following is a version of \cite[Proposition 3.3]{BS2}.

 
 \bP\label{BS2Prop3p3Nov25}
 Suppose $A$ is a noetherian commutative ring, with $\dim A=d$ 
 and $P$ is a projective $A$-module, with $rank(P)=n$.
 Assume $2n\geq d+3$ and $P\cong Q\oplus A$. 
Let $J_0, J_1, J_2, J_3\subseteq A$ be  ideals, with $height(J_i)\geq n$ for $i=0, 1, 2, 3$, 
$J_0+J_1J_2=A$ and 
$J_0J_1J_2+J_3=A$.
Also, let
$$
\alpha: P \sur J_0\cap J_1, \quad \beta: P\sur J_0\cap J_2\quad {\rm be ~surjective~maps}~\ni \alpha\otimes \frac{A}{J_0} = \beta\otimes \frac{A}{J_0}.
$$
Further, assume that 
 there is a surjective map 
$$
\gamma: P \sur J_1\cap J_3 \quad \ni \quad \gamma\otimes \frac{A}{J_1}=\alpha\otimes \frac{A}{J_1}.
$$
Then, there is a surjective map
$$
\delta: P \sur J_2\cap J_3 \quad \ni \quad \delta \otimes \frac{A}{J_3}= \gamma \otimes \frac{A}{J_3}, \quad \delta \otimes \frac{A}{J_2}= \beta \otimes \frac{A}{J_2}.
$$

If $A=R[X]$ is a polynomial ring over a regular ring $R$, over an infinite field $k$, same is true, when $2n\geq \dim A+2$.
 \eP
 \pf Denote $\omega_0=\alpha\otimes \frac{A}{J_0}=\beta\otimes \frac{A}{J_0}$, $\omega_1=\alpha\otimes \frac{A}{J_1}=\gamma\otimes \frac{A}{J_1}$
 $\omega_2=\beta\otimes \frac{A}{J_2}$, $\omega_3=\gamma\otimes \frac{A}{J_3}$.
 By Moving Lemma \ref{movingLemm} there is ${\bf u}=(f, p, s)\in \widetilde{{\CQ}}(P)$, such that 
 $$
 \eta({\bf u})= (J_0, \omega_0),\quad \eta(\Gamma({\bf u}))=(J_4, \omega_4), \quad J_1J_2J_3+J_4=A, \quad height(J_4)\geq n.
 $$
 %
 As is intended, 
 $f(P)=J_0\cap J_4$, with $J_0+J_4=A$. 
 
Denote 
$g:=f:P\sur J_0J_4$ be the 
 a surjective map defined by $f$. 
 Then, $g\otimes \frac{A}{J_0}=\omega_0$.
 By Addition Principle \cite[Theorems 5.6, 5.7]{BK}, applied to $\gamma$ and $g$, there is a surjective map
 $$
 \mu: P \sur (J_1\cap J_3)\cap (J_0\cap J_4) ~~ \ni \quad \mu \otimes \frac{A}{J_1\cap J_3}=\gamma\otimes \frac{A}{J_1\cap J_3}= \omega_1\star\omega_3,
 $$
 $$
 {\rm and} \qquad \mu \otimes \frac{A}{J_0\cap J_4}=g \otimes \frac{A}{J_0\cap J_4}=  \omega_0\star\omega_4.
 $$
 It follows, $\mu \otimes \frac{A}{J_0\cap J_1}=\omega_0\star \omega_1= \alpha\otimes \frac{A}{J_0\cap J_1}$. 
 By Subtraction Principle  \cite[Theorems 3.7, 4.11]{BK}, applied to $\mu$ and $\alpha$, there is a surjective map 
$\nu: P\sur J_3\cap J_4$
 such that $\nu\otimes \frac{A}{J_3\cap J_4} = \mu \otimes \frac{A}{J_3\cap J_4}=\omega_3\star \omega_4$.
 By Addition Principle \cite[Theorems 5.6, 5.7]{BK}, applied to $\nu$ and $\beta$, there is a surjective map 
 $$
 \lambda:P \sur (J_0\cap J_2) \cap (J_3\cap J_4)
 \quad \ni \quad \lambda\otimes \frac{A}{J_0\cap J_2}= \beta \otimes \frac{A}{J_0\cap J_2}=\omega_0\star\omega_2
 $$
 $$
 {\rm and,}\quad ~~\lambda\otimes \frac{A}{J_3\cap J_4}= \nu \otimes \frac{A}{J_3\cap J_4}= \omega_3\star \omega_4.
$$
Now apply Subtraction Principle \cite[Theorems 3.7, 4.11]{BK}, to $\lambda$ and $g$. There is a surjective map
$$
\delta: P\sur J_2\cap J_3 \quad \ni \quad \delta\otimes \frac{A}{J_2\cap J_3} = \lambda \otimes \frac{A}{J_2\cap J_3}=\omega_2\star \omega_3.
$$
So, $\delta\otimes \frac{A}{J_2} =\omega_2$ and  $\delta\otimes \frac{A}{J_3} =\omega_3$.
\pic $\eop$

 The following is the version of Corollary \ref{newEuler}.
 \bT\label{oldGoldMethod}
 Suppose $A$ is a commutative noetherian ring with $\dim A=d$ and $P$ is a projective $A$-module, with $rank(P)=n$.
 Assume $2n\geq d+3$ and $P\cong Q\oplus A$. 
 Let $(J, \omega_J) \in {\mathcal LO}^n(P)$ and $\overline{\left(J, \omega_J\right)}=0\in E(P)$.
 Then, $\omega_J$ lifts to a surjective map $P\sur J$.
 
If $A=R[X]$ is a polynomial ring over a regular ring $R$, over an infinite field $k$, same is true when $2n\geq \dim A+2$.

 \eT
 

\pf Suppose 
$(J, \omega_J) \in {\mathcal LO}^n(P)$ and $\overline{\left(J, \omega_J\right)}=0\in E(P)$.
We have a set 
$$
\left\{(J_t, \omega_t) : 1 \leq t \leq r+s \right\} 
$$
such that 
\begin{enumerate} 
\item $height (J_t)=n.$ 
\item there are surjective maps 
$\alpha_t:  P \sur J_t$ such that $\alpha_t$ lifts $\omega_t.$ 

\item And 
\begin{equation}\label{eulerEqnNov25}
(J, \omega)_{{\BZ}}+\sum_{l=r+1}^{r+s} (J_t, \omega_t)_{{\BZ}} = \sum_{t=1}^r  (J_t, \omega_t)_{{\BZ}} 
\quad {\rm in}\quad \BZ\left({\mathcal LO}^n_c(P)\right).
\end{equation}
holds in the free group $ \BZ\left({\mathcal LO}^n_c(P)\right)$.
\end{enumerate}
First assume that $J_1, J_2, \ldots, J_r$ are pairwise comaximal.
In this case, $J, J_{r+1}, \ldots, J_{r+s}$ are pairwise comaximal. 
Write 
$$
J'=\cap_{l=r+1}^{r+s} J_t, \qquad J"=\cap_{t=1}^r J_t.  
\qquad {\rm Then}\quad J\cap J' = J".
$$  
Further, by Addition Principle \cite[Theorems 5.6, 5.7]{BK}, there are surjective homomorphisms 
$\alpha': P\sur J'$  and $\alpha": P \sur J"$ such that
$$
(J', \omega')_{{\BZ}} = \sum_{t=r+1}^{r+s}(J_t, \omega_t)_{{\BZ}},  \quad 
(J", \omega")_{{\BZ}} = \sum_{t=r+1}^{r+s}(J_t, \omega_t)_{{\BZ}} 
\quad {\rm in}\quad \BZ\left({\mathcal LO}^n_c(P)\right)
$$
where $\omega' =\alpha'\otimes A/J'$ and $\omega" =\alpha"\otimes A/J"$.
 So, by Subtraction Principle \cite[Theorems 3.7, 4.11]{BK},
there is a surjective homomorphism $\alpha: P\sur J$ such that 
$\alpha \otimes A/J =\alpha"\otimes A/J= \omega.$

Now, we consider that $J_1, J_2, \ldots, J_r$ are not, necessarily, pairwise 
comaximal. Given an Equation, as in (\ref{eulerEqnNov25}), we would associate an integer $n({\rm Eqn-}\ref{eulerEqnNov25})\geq 0$, as follows.
Let $S_i$ be the set of all connected components of $J_i$ and 
$S=\cup_{i=1}^rS_i.$ 
For $K \in S,$ let $n(K)+1$ be the cardinality
of the set $\{t: K+J_t \neq A \}.$ Let $n({\rm Eqn-}\ref{eulerEqnNov25})= \sum_{K\in S} n(K).$
We have $n({\rm Eqn-}\ref{eulerEqnNov25})=0$ if and only if $J_1, J_2, \ldots, J_r$ are comaximal. 
 
Now,  assume $n({\rm Eqn-}\ref{eulerEqnNov25})\geq 1$. Therefor, $n(K)\geq 1$ for some $K \in S.$ 
We can assume $K \in S_1$  and $K+J_2 \neq A.$ 
So, $\exists ~ \tilde{K}$ a connected component of $J_2$ such that 
$K+ \tilde{K} \neq A.$    

First, assume $K \neq \tilde{K}.$ Both $K,\tilde{K}$ cannot be connected component of $J.$ 
 ({\it components add up to $A.$})
Without loss of generality, assume $K$ is not a
connected component of $J.$ 
Using Eqn-\ref{eulerEqnNov25}, it follows that there is an integer $l$, with $r+1 \leq l \leq r+s$, 
 such that (1) $K$ is a connected component of $J_l$, (2) $\alpha_l \otimes A/K = \alpha_1 \otimes A/K.$  
Assume $l=r+1$ and denote $\omega_K:= \alpha_l \otimes A/K = \alpha_1 \otimes A/K: P \sur K/K^2$.
%
%
We  write $J_1=K \cap K_1$ and $J_{r+1}= K \cap K_2$
where $K+K_1=A=K+K_2.$ By Moving Lemma \ref{movingLemm}, applied to $\omega_{K_1}:=\alpha_1\otimes \frac{A}{K_1}$,
there is an ideal $K_3$ such that (3) $height(K_3)\geq n$, (4) $K_3$ is comaximal to $J, J_j, \forall 1 \leq j \leq r+s$,
(5) there is a surjective map  $\beta: P \sur K_3 \cap K_1$  such that $\alpha_1 \otimes A/K_1= \beta \otimes A/K_1.$
\\
We have three surjective maps:
$$
\alpha_1: P \sur K \cap K_1, \quad \alpha_{r+1}: P \sur K \cap K_2 \quad \beta: P \sur K_1\cap K_3
$$
By proposition \ref{BS2Prop3p3Nov25}, there is a surjective map 
$$
\beta_{r+1}: P \sur K_3 \cap K_2\quad \ni ~~
\alpha_{r+1}\otimes \frac{A}{K_2}=\beta_{r+1}\otimes \frac{A}{K_2}, 
 ~~\beta\otimes \frac{A}{K_3} 
=\beta_{r+1}\otimes\frac{A}{K_3}. 
$$ 
So, we have 
\begin{equation}\label{EupfEquanTwo}
(J, \omega)_{{\BZ}}+\left(\widetilde{J_{r+1}}, \widetilde{\beta_{r+1}} \right)_{{\BZ}}
+\sum_{l=r+2}^{r+s}(J_l, \omega_l)_{{\BZ}} 
 =
\left(\widetilde{J_{1}}, \widetilde{\beta_{1}} \right)_{{\BZ}} + 
\sum_{l=2}^{r}(J_l, \omega_l)_{{\BZ}}.  
\end{equation}
where $\widetilde{J_{r+1}}= K_3 \cap K_2$ and $\widetilde{J_{1}}= K_3 \cap K_1.$   
({\it $K$ is removed from both sides and $K_3$ is inserted.})
It is clear $n({\rm Eqn-}\ref{EupfEquanTwo})<n({\rm Eqn-}\ref{eulerEqnNov25}).$
Therefore, by induction, the  Equation-\ref{EupfEquanTwo} would reduce to an Equation (*), so that 
$n(*)=0$.

Now assume $K= \tilde{K}.$ Let $\omega_K= \alpha_1 \otimes A/K.$ 
Therefore, $K=\tilde{K}$ is a component of $J_2$. We denote $\tilde{\omega}_K
=\alpha_2 \otimes A/K.$   Using Equation-\ref{eulerEqnNov25}, it follows that either  $(K, \omega_K)$ or $(K, \tilde{\omega}_K)$ 
is a summand of $\sum_{t={r+1}}^{r+s}(J_t, \omega_t)_{{\BZ}}$. Without loss of generality, we assume that 
$(K, \omega_K)$ is a summand of $(J_{r+1}, \omega_{r+1})_{{\BZ}}$ and complete the induction exactly in the same manner, as above.
\tcp
$\eop$


\subsection{Comparison with Chow Groups} \label{compChowSec}
In this section, we exploit the work of N. Mohan Kumar and M. P. Murthy  \cite{MoM, Mu}
to compare  the Euler class group $E(P)$ with the Chow group 
of zero cycles $CH^d(A)$, when $A$ is a smooth affine algebra over an algebraically closed field and $n=rank(P)=d$.

\bD\label{euler2F0k0}{\rm 
Let $A$ be a Cohen Macaulay ring, with $\dim A=d$. Let $K_0(A)$ denote the Grothendieck group of projective $A$-modules. Let
$F_0K_0(A)=$
$$
\left\{ \left[\frac{A}{I} \right]\in K_0(A): I~{\rm is ~a~ locally~ complete~ intersection~ ideal, ~with}~
height(I)=d
\right\}.
$$
It was established in \cite[Theorem 1.1]{M3} that $F_0K_0(A)$ is a subgroup of $K_0(A)$. 

Let $Q$ be a projective $A$-module with $rank(Q)=d-1$, and $P=Q\oplus A$. Then, there is a surjective homomorphism 
$$
\varphi: E(P) \lra F_0K_0(A)\quad {\rm sending} \quad [(I, \omega)] \mapsto \left[\frac{A}{I}\right]
$$
}
\eD
\pf Since $A$ is Cohen Macaulay, for $(I, \omega) \in {\mathcal LO}^d(P)$, $I$ is
a  locally complete intersection ideal. Now, consider the map ${\mathcal LO}^d_c(P) \lra F_0K_0(A)$, sending $(I, \omega) 
\mapsto \left[\frac{A}{I}\right]$. This map extends to a homomorphism $\varphi_0:{\BZ}\left({\mathcal LO}^d_c(P)\right) \lra F_0K_0(A)$. 
Now, if $(I, \omega)$ is a global orientation, then $\omega$ lifts to a 
surjective map $f:P\sur I$. Since $P=Q\oplus A$, it follows 
$$
\varphi_0(I, \omega)= \left[\frac{A}{I}\right]=\sum_{r=0}^d(-1)^r\left[\wedge^rP \right]=0.
$$
Therefore, $\varphi_0$ factors through a map $\varphi: E(P) \lra F_0K_0(A)$.
For, $ \left[\frac{A}{I}\right]\in F_0K_0(A)$, there is an isomorphism $\frac{P}{IP} \iso \frac{I}{I^2}$,
which gives rise to a surjective map $\omega: P \sur \frac{I}{I^2}$. Therefore, 
$\varphi([(I, \omega)])=
\left[\frac{A}{I}\right]$. So, $\varphi$ is surjective.

\vspace{2mm}
Now, we assume that the base field $k$ is algebraically closed, and use \cite{MMu}.
\bC\label{CorSwanBertini}
Suppose $A$ is a reduced affine algebra over an algebraically closed field $k$, with $\dim A=d\geq 2$. Assume $A$  is Cohen Macaulay and that $F_0K_0(A)$ has no $(d-1)!$ torsion.
Let $Q$ be a projective $A$-module with $rank(Q)=d-1$, and $P=Q\oplus A$.
Then, the map $\varphi: E(P) \lra F_0K_0(A)$ in (\ref{euler2F0k0}) is  an isomorphism.
\eC
\pf By Swan's Bertini theorem \cite[pp. 586]{MoM}, it follows that $F_0K_0(A)$ coincides with the usual subgroup $F^dK_0(A)$ (see \cite{F, Mu, MMu}), 
which is generated by the cycles of $A/{\m}$, where $\m$ runs through the smooth maximal ideals
of height $d$. We only need to prove that $\varphi$ is injective. Suppose $\varphi(x)=0$ for 
some $x\in E(P)$. By Moving Lemma, we can write $x=[(I,\omega)]$, for some $(I,\omega)\in {\mathcal LO}^d(A)$. 
Therefore 
$\varphi(x)=\left[\frac{A}{I} \right]=0$. Since $P=Q\oplus A$, 
the top Chern class $C^d(P^*)= \sum_{r=0}^d (-1)^r[\wedge^r P]=0\in F^dK_0(A)$  (see \cite[Definition 3.5]{Mu}).
Therefore, $C^d(P^*)=  \left[\frac{A}{I} \right]\in F^dK_0(A)$.
By  \cite[Theorem 2.1]{MMu}, it follows $\omega$ lifts to surjective map $P\sur I$. Therefore, $x=[(I, \omega)]=0$. This establishes that $\varphi$ is an isomorphism. \pic $\eop$


\vspace{2mm}
The condition in (\ref{CorSwanBertini}) that $F_0K_0(A)$ has no $(d-1)!$ torsion is a minor condition, due to the results of Levine \cite{Le} and Srinivas \cite{S}
(see  \cite[Lemma 2.10, Theorem 2.14]{Mu}).
Summarizing all the above, with smoothness hypotheses, we have the following.
\bC\label{CorSallFour}
Suppose $A$ is smooth affine algebra over an algebraically closed field $k$, with $1/2\in k$ and $\dim A=d\geq 3$. 
Let $Q$ be a projective $A$-module with $rank(Q)=d-1$, and $P=Q\oplus A$.
Then, the maps 
$$
\diagram
\pi_0\left({\mathcal LO}(P)\right)
&
E(P) \ar[l]_{\qquad\sim}\ar[r]^{\sim\quad}
&
F_0K_0(A) &
 CH^d(A) \ar[l]_{\sim}\\
\enddiagram
$$
are isomorphisms, where $CH^d(A)$ denotes the Chow group of codimension $d$ cycles. 
\eC
\pf The last isomorphism  follows from Riemann-Roch theorem, because
 $F_0K_0(A)$ is divisible and  does not have $(d-1)!$ torsion \cite[Lemma 2.10, Theorem 2.14]{Mu}. 
The second isomorphism follows from (\ref{CorSwanBertini}), while the first isomorphism follows from (\ref{useBKNov14}).
\pic $\eop$

\subsection{Some Closing Remarks}\label{closingSec}
Before we close this main body of this article, we have the following remarks. 

\begin{remark}{\rm

For the following comments,  
assume $A$ is an essentially smooth affine rings, over an infinite perfect field $k$, with $1/2\in k$ and $\dim A=d\geq 3$, and $X=\spec{A}$.
%
\bE
\item The structure theorem \cite[Theorem 4.21]{BDM} illustrates that these 
monoids $\pi_0\left(\widetilde{{\CQ}}(P)\right)$ can assume a wide range of values.
%

\item\label{loose} Assume $P$ does not have a unimodular element (see \cite{Mk}). 
Then, there is no ideal preserving and homotopy preserving map ${\mathcal LO}(A^n) \lra {\mathcal LO}(P)$.
%
\bE
\item\label{allGroupstopRank} However, in a subsequent article \cite{MM2}, we prove that when $rank(P)=d$, then $\pi_0\left({\mathcal LO}(P)\right) \cong 
\pi_0\left({\mathcal LO}(\Lambda^dP \oplus A^{d-1})\right)$. Since the latter one is a group, $\pi_0\left({\mathcal LO}(P)\right)$ is a group. 
In particular, for $\pi_0\left({\mathcal LO}(P)\right)$ to be a group, it is not necessary that $P$ splits as $P\cong Q\oplus A$.
\item {\bf Open Problem:} It remains open, whether $\pi_0\left({\mathcal LO}(P)\right)$ is always a group or not, whenever $2n\geq d+3$.


\item In \cite{MM2}, we establish a natural additive map $\pi_0\left({\mathcal LO}(P)\right) \lra CH^n(A)$, where $CH^n(A)$ denotes the Chow group 
of codimension $n$ cycles.
\eE
\item If $\pi_0\left({\mathcal LO}(P) \right)$ is a group and ${\bf e}_0\neq {\bf e}_1$,
  then the natural map ${\mathcal LO}(P)\lra \pi_0\left({\mathcal LO}(P) \right)$ (see (\ref{EntopizeroNov14})),
 does not factor through
a group homomorphism, from $E(P)$ to $\pi_0\left({\mathcal LO}(P) \right)$. This is because the global orientations, map to 
${\bf e}_0$.

\item Note (see \cite{F}), that the total Chern class $C(P)=1+C^1(P)+\cdots+C^n(P)$ takes value the total Chow groups 
$CH(X)=\oplus_{i=1}^dCH^i(X)$, which is an invariant of $X$. However, the homotopy obstruction group 
$\pi_0\left({\mathcal LO}(P)\right)$, which houses the Nori class $\varepsilon(P)$, is an invariant of $P$. 
This does not come as a surprise, because when $rank(P)=n=d$, the Euler class of $P$ (as in \cite{BS3}) 
takes value in $E^n(A, \wedge^dP)=E(\wedge^dP\oplus A^{d-1})$, which is dependent on $P$.
\eE
}
\end{remark} 

\appendix
\section{The  Motivic Interpretation}\label{motivicSec}
In this section, we attempt to give a motivic interpretation to the homotopy obstruction sets, in analogy to the  case 
when $P=A^n$  \cite{BM, Mo}.
Four descriptions for the same was given in section \ref{foundationSec}, assuming $1/2\in A$. 
For our purpose, in this section, it would be best to work with $\widetilde{{\CQ}}'(P)$ and $\pi_0\left(\widetilde{{\CQ}}'(P)\right)$. 
We assume $1/2\in A$ in this section. Recall,
with
\begin{equation}\label{motivB2n}
 \SB_{2n+1}=\frac{k[X_1, \ldots, X_n; Y_1, \ldots, Y_n, Z]}{\left(\sum_{i=1}^nX_iY_i+Z^2-1\right)}, \quad{\rm and}\quad
 Q_{2n}'=\spec{\SB_{2n+1}}
\end{equation}
$[\spec{A}, Q_{2n}']_{\Sch}\cong$
$$
Q_{2n}'(A)\cong \left\{(f_1, \ldots, f_n; g_1, \ldots, g_n, z)\in A^{2n+1}:
\sum_{i=1}^nf_ig_i+z^2=1 \right\}.
$$
Also recall, $\pi_0\left(Q_{2n}' \right)(A) \cong [Q_{2n}', \spec{A}]_{{\BA}^1}$ where the right hand side 
denotes the set of  all morphisms in  the
${\BA}^1$-homotopy 
category \cite[Chapter 8]{Mo} (also see \cite[Theorem 1.1.1]{AF}). 
A similar interpretation for $\widetilde{{\CQ}}'(P)$ and $\pi_0\left(\widetilde{{\CQ}}'(P)\right)$ would be desirable.



We follow Swan \cite[\S 1, 2]{Sw}. Suppose $Q$ is a projective $A$-module. Let
$S(Q^*)=\bigoplus_{i\geq 0} S_i(Q^*)$ denote the symmetric algebra 
of $Q^*$. 
Let $Quad(Q)=\{\varphi\in Hom(Q, Q^*): \varphi^*=\varphi \}$ denote the $A$-module of all the quadratic forms on $Q$.
Given $\varphi\in Quad(Q)$, let $B(\varphi) \in Hom(Q \otimes Q, A) \cong Q^*\otimes Q^*$ be the corresponding bilinear map.
In fact, this association $\varphi \mapsto B(\varphi)$ induces a 
bijection $Quad(Q)\iso S_2(Q^*)$ (see  \cite[\S ~2]{Sw}).

Since $A$ is commutative, all maps $f:Q^*\lra A$ extends to a map $S(Q^*) \lra A$. So, we have the commutative diagram of 
bijections:
 $$
 \diagram
 Q \ar[r]^{ev\qquad}_{\sim\qquad}\ar[dr]_{\lambda}& Hom(Q^*, A)\ar[d]^{\wr} \\
 &  Hom(S(Q^*), A) \\
 \enddiagram
 $$
  $$
  {\rm For}~x\in Q,~~
  f, g\in Q^*\quad \langle \lambda(x), f\rangle =f(x)\quad \langle \lambda(x), fg\rangle =f(x)g(x)
  $$ 
For a bilinear map $\beta\in Hom(Q\otimes Q, A)=Q^*\otimes Q^*$, we can write $\beta=\sum f_i\otimes g_i$ for some $f_i, g_i\in Q^*$.
and
   $$
   \langle \lambda (x), \beta\rangle =\sum f_i(x)g_i(x)=\beta(x, x).
   $$

  Fix a quadratic form $\varphi: Q\lra Q^*$ and $B(\varphi): Q\otimes Q\lra A$ 
  be the corresponding bilinear map. More precisely, $B(\varphi)(x, y)=\varphi(x)(y)$. As usual, define $q: Q \lra A$ by $q(x)=B(x, x)$.
   Then,
   $$
   {\rm for}~x\in Q\quad   \langle \lambda(x), B(\varphi)\rangle =B(\varphi)(x, x)=q(x).
   $$
  We introduce some notations.
   \begin{notations}\label{notaBqSwan}{\rm 
    Suppose $A$ is a commutative noetherian ring, 
    with $1/2\in A$, and $X=\spec{A}$.
For a  quadratic space $(Q, \varphi)$ over $A$, denote
   $$
   \BS(\varphi)=\{x\in Q: q(x)=1 \}, ~
  \SB(\varphi)=\frac{S(Q^*)}{( B(\varphi)-1)}
  ~{\rm and} ~ \CX(\varphi)=\spec{\SB(\varphi)}.
   $$
   }
   \end{notations}
   
   \bP\label{swanSphere}
  With notations as in (\ref{notaBqSwan}), the following maps
   $$
 \diagram
  [X,{\CX}(\varphi)]_{\Sch_A}  \ar[r]^{\sim} &Hom\left(\SB(\varphi), A\right) &  \BS(\varphi) \ar[l]_{\qquad\quad\sim}\\
  \enddiagram
    \quad {\rm are~bijections},
   $$
   where $[-, -]_{\Sch_A}$ denotes the set of morphisms in $\Sch_A$.
   \eP
   \pf Follows from above discussions. $\eop$
   
   \vspace{3mm}
   \begin{remark}{\rm 
   Use the same notations, as in (\ref{notaBqSwan}). Consider the pre sehaf
   $$
[-,{\CX}(\varphi)]_{\Sch}:  \underline{Sch}_{A}\lra \underline{Sets} \quad {\rm sending}\quad 
  Y \mapsto [Y,{\CX}(\varphi)]_{\Sch_A} 
   $$
   In fact, for affine schemes $Y=\spec{B}\in \underline{Sch}_{A}$, the following maps %
   $$
\diagram
[Y,{\CX}(\varphi)]_{\Sch_A}\ar[r]^{\sim} &   Hom\left({\SB}(\varphi),  B\right) 
    &   \BS\left(\varphi\otimes B\right)\ar[l]_{\qquad\sim}\\
\enddiagram
\quad {\rm   are~ bijections.}
   $$

   One can make a similar statement for any scheme $Y\in {\Sch}_A$.  Let $f: Y \lra X$ be the structure map, and $f^*$ would denote the pullback. Redefine
   $$
   \left\{
   \begin{array}{l}
   {\BS}(f^*q)=\{x\in \Gamma(Y, f^*Q): f^*q(x)=B(f^*\varphi)(x, x)=1 \}\\
     \SB(f^*\varphi)=\frac{S(f^*Q^*)}{( B(f^*\varphi)-1){\CO}_Y}\\
 \CX(f^*\varphi)=\spec{\SB(f^*\varphi)}.\\
   \end{array}
   \right.
   $$
   Then,
   the following maps %
   $$
  \diagram
 [Y,{\CX}(\varphi)]_{\Sch_A}  \ar[r]^{\sim\qquad} &Hom\left({\SB}(\varphi),   \Gamma(Y, {\CO}_Y)\right)&      \BS\left(f^*\varphi\right)\ar[l]_{\qquad\quad\sim} \\ 
\enddiagram
   $$
are bijections. (see \cite[II, Ex 2.4]{H}).
   }
   \end{remark}
   \vspace{3mm}
   Corresponding to the notations (\ref{motivB2n}), we introduce the following notations. 
   \begin{notations}\label{notaMotivBP}
   Let $({\BQ}(P), q) = {\BH}(P)\perp A$ be as in (\ref{nota}). Denote the underlying projective module of $({\BQ}(P), q)$ 
   by the same notaion ${\BQ}(P):=P^*\oplus P \oplus A$.
   Let $B: {\BQ}(P)\times {\BQ}(P) \lra A$ be the corresponding bilinear form.
   Define
    $$
   \SB(P)=\frac{S(P\oplus P^* \oplus A)}{(B-1)}, \quad {\rm and~denote}\quad Q_P':=\spec{\SB(P)}.
   $$
   
   \end{notations}
   
   \bC\label{ourSphere}
   %
   %
   Use the notations, as in (\ref{notaMotivBP}).
   Then, for $Y\in \Sch_A$ and the structure map $f:Y \lra X$, the following maps %
   $$
   \left\{
   \begin{array}{l}
    \widetilde{\CQ}'(P)={\BS}(q) \iso
    \left[X, Q_P'\right]_{\Sch_A}\\
  %
 \widetilde{\CQ}'(f^*P)={\BS}(f^*q) \iso 
     \left[Y, Q_P'\right]_{\Sch_A} \\
    \end{array}
    \right.
    \quad {\rm are~bijections}.
   $$
   Consequently, the association 
   $$
   Y \mapsto  \widetilde{\CQ}'(f^*P)\quad{\rm defines ~a~ pre sheaf} \quad\Sch_A\lra \Sets.
   $$
   Therefore,
   one can define 
   $$
   \pi_0(\widetilde{\CQ}(P)): \Sch_A \lra \Sets  \quad {\rm as ~a~pre sheaf}.
   $$
  
   \eC
 \pf Follows from (\ref{swanSphere}).
 \tcp $\eop$.


\vspace{2mm}
In analogy to to the free case $P=A^n$, we raise the following question.

\begin{question}\label{motivicQn}
{\rm 
Suppose $k$ is an infinite perfect field, with $1/2\in k$ and $A$ is an essentially smooth ring over $k$, with $\dim A=d$. Write $X=\spec{A}$.
Suppose
$P$ is a projective $A$-module with $rank(P)=n$.
%
The question remains, whether a motivic interpretation can be given to the pre sheaf $ \pi_0(\widetilde{\CQ}(P))$.
In particular, whether 
$$
  \pi_0(\widetilde{\CQ}(P))(A) \cong  \left[X, Q_P\right]_{{\BA}^1},
$$
where $\left[X, Q_P\right]_{{\BA}^1}$ denotes the set of  all maps $X\lra Q_P$, in the $\BA^1$-homotopy category?
It may be best to assume $n\gg 0$.

If $A$ is local, then  the  equality holds (see \cite[Chapter 8]{Mo}).
}
\end{question}
\end{document}